\documentclass{amsart}
\usepackage[utf8]{inputenc}
\usepackage{amssymb,amsmath,amsthm}
\usepackage[textwidth=13cm,textheight=21cm]{geometry}
\usepackage{xcolor}
\usepackage[colorlinks=true,linkcolor=blue]{hyperref}
\usepackage{siunitx}
\usepackage{tikz-cd}
\usepackage[shortlabels]{enumitem}
\usepackage{mathtools}

\newtheorem{thm}{Theorem}[section]
\newtheorem{lem}[thm]{Lemma}
\newtheorem{defin}[thm]{Definition}
\newtheorem{prop}[thm]{Proposition}
\newtheorem{cor}[thm]{Corollary}
\newtheorem{rem}[thm]{Remark}

\newcommand{\N}{\mathcal{N}}

\newcommand{\M}{{\mathcal{M}}}

\newcommand{\hol}{\mathrm{hol}}
\newcommand{\sm}{\smallsetminus}
\newcommand{\dev}{\mathrm{dev}}
\newcommand{\HH}{{\mathcal{H}}}
\newcommand{\R}{{\mathbb{R}}}

\newcommand{\C}{{\mathbb{C}}}
\newcommand{\Z}{{\mathbb{Z}}}
\newcommand{\Mod}{{\operatorname{Mod}}}
\newcommand{\PMod}{{\operatorname{PMod}}}

\newcommand{\q}{\boldsymbol{q}}
\newcommand{\df}{\overset{\operatorname{def}}{=}}
\newcommand{\HHm}{\mathcal{H}_{\mathrm{m}}}

\newcommand{\GL}{\mathrm{GL}_2^+(\R)}
\newcommand{\SL}{\mathrm{SL}_2(\R)}
\newcommand{\U}{\mathcal{U}}

\newcommand{\W}{\mathcal{W}}

\title[Geometric structures, orbit closures and invariant measures]{Geometric structures, orbit closures and invariant measures in moduli spaces}

\author{John Smillie}
\address{University of Warwick 
{\tt j.smillie@warwick.ac.uk}}
\author{Florent Ygouf}
\address{Université de Rennes
{\tt florent.ygouf@univ-rennes.fr}}

\begin{document}

\begin{abstract}Building on the seminal work of Eskin,  Mirzakhani and Mohammadi in the context of dynamics on moduli spaces, we introduce linear loci and linear measures to provide a relatively basic description of orbit closures and invariant measures in strata of translation surfaces. We use this to discuss geometrical properties of orbit closures and invariant measures from the perspective of geometric structures. 
\end{abstract}

\maketitle

\tableofcontents

\newpage

\section{Introduction}
A fundamental insight in the study of translation surfaces is that many dynamical and geometric questions can be reduced to understanding the $\GL$-action on strata of translation surfaces. Some questions involve understanding the horocycle flow on $\GL$-orbit closures. In \cite{EM,EMM} Eskin, Mirzakhani and Mohammadi study the $\GL$-action on finite manifold covers of strata and obtain strong structural restrictions on orbit closures and invariant measures.  In this paper, we build on their work to give two distinct descriptions of orbit closures and invariant measures in strata. 

The need for the first description arose when writing \cite{SSWY}. This paper deals with horospherical measures, extending earlier work of \cite{LM} which dealt with the principal stratum. Horospherical measures can be described loosely as horocycle invariant measures with maximal additional invariance. An alternative characterization of the structure of orbit closures was formulated in \cite[Theorem 2.10]{SSWY} and does not exactly match that of \cite{EM,EMM}. The proof of this result was deferred to the present paper where it appears as Theorem \ref{thm: orbit is invariant subvariety}, which is itself a consequence of Theorem \ref{thm: orbit is linear}. Theorem \ref{thm: invariant measure is cone measure associated to linear measure} establishes an analogous characterization for invariant measures. These results rely on the notions of {\em linear loci} and {\em linear measures} which we introduce here. 

To explain the purpose of our second description of orbit closures and invariant measures, some context is needed. According to \cite{EMM}, $\GL$-orbit closures in the unnormalized stratum $\HH$ behave very much like submanifolds of strata which are locally described by sets of linear equations with respect to the period coordinate charts. However, they can fail to be manifolds in two distinct ways: they can have orbifold points and they can have self-intersections. (See example 4.1.10 in \cite{Filip_survey} for an example of an orbit closure with self-intersections.) Eskin, Mirzakhani and Mohammadi deal with the first problem by passing to a finite cover of the stratum $\widetilde\HH$ of $\HH$ which has no orbifold points and restricting their attention to orbit closures of the $\GL$-action on this finite cover. They deal with the second problem by ``desingularizing'' orbit closures \textit{i.e.} by showing that each orbit closure is the image of a manifold by a map $f:\N \to\widetilde\HH$ which is generically injective, proper and ``locally linear" in the  stratum $\HH$. However, in both \cite{EM} and \cite{EMM} the theorems are {\em stated} as if they applied to strata rather than to finite covers. The fact that the proofs require taking a finite cover of the stratum is mentioned in \cite{EM} but is not mentioned in \cite{EMM}. Thus the results of \cite{EM} and \cite{EMM} only apply to finite manifold covers of strata and not to strata themselves. (See Section \ref{sec: recasting} for further discussion). The purpose of our second description is threefold. It addresses this finite cover issue, endows the map $f$ with useful additional structure and makes more explicit the connection of the work of Eskin, Mirzakhani and Mohammadi with the language of affine structures, a connection already highlighted by Quint in \cite{Quint}. We show that $\N$, like the stratum $\HH$, carries an affine structure. In fact, we show that this structure is richer than merely affine: it has a linear holonomy compatible with a real structure. We capture this extra feature with our notion of ``real linear orbifolds''. 

\begin{thm}\label{thm: linear structure on orbit closures}
    Let $\M\subset \HH$ be a $\GL$-orbit closure. There is a connected real linear orbifold $\N$ endowed with an action of $\GL$ and a proper equivariant real linear immersion $f : \N \to \HH$ whose image is $\M$ and which is injective on an open dense set.  
\end{thm}

One useful feature of this description is the fact that $\N$ is endowed with a $\GL$-action for which the map $f$ is equivariant. 
This suggests that one could approach dynamical problems by studying them in the ``desingularization" $\N$ where issues of self-intersection do not arise.  Such an approach is used in an ongoing project of Arana-Herrera, De Witt, Eskin, Gadre, Guttirez-Romo, Lima, Matheus, Rafi and Schleimer. They study the algebraic hull of the Kontsevich-Zorich cocycle on orbit closures via the flow group and its associated monodromy.  The orbifold fundamental group of $\N$ provides a natural candidate for a fundamental group associated to the orbit closure. 

The characterization corresponding to Theorem \ref{thm: linear structure on orbit closures} for invariant measures using the language of geometric structures relies on the notion of real equilinear orbifolds. These are complex affine orbifolds which come with a natural family of Radon measures that can be understood as a Lebesgue measure in charts. We refer to Subsection \ref{sec: linear structures} for precise definitions. The cone construction (\cite{Athreya-Masur} Subsec. 3.4) gives a method of inducing a measure on $\HH^1$ from a measure on $\HH$. We will give a generalization of the cone construction and relate it to the approach of \cite{EM} in Subsection \ref{sec: from H to H1}.

\begin{thm}\label{thm: linear structure on invariant measures} 
    Let $\nu$ be an ergodic $\SL$-invariant probability measure on $\HH^1$. There is a connected real equilinear orbifold $\N$ endowed with an action of $\GL$ by real linear maps, an ergodic homogeneous Lebesgue measure $\lambda$ on $\N$ and a proper  equivariant real linear immersion $f : \N \to \HH$ such that $\nu$ is obtained from $f_*\lambda$ by the cone construction.
\end{thm}

\subsection{Organization of the paper}

In Section \ref{sec: prelims}, we introduce the main objects needed in this paper and set up our notation. We also recall the interplay between the action of $\GL$ on $\HH$ and that of $\SL$ on $\HH^1$. In Section \ref{sec: recasting} we carefully describe the results of \cite{EM} and \cite{EMM} that we are using as the basis for our results.  We introduce our notion of locally linear sets and locally linear measures and use them as a vehicle for stating the part of the results of \cite{EM} and \cite{EMM} that we need. In Section \ref{sec: linear loci and linear measures}, we introduce our notion of linear loci and linear measures and show that a locally linear set is a linear locus and a locally linear measure is a linear measure. We then use these results and the results of Section \ref{sec: recasting} to show that that $\GL$-orbit closures in $\HH$ are irreducible linear loci and $\SL$-invariant ergodic measures in $\HH^1$ are cone measures associated to irreducible linear measures. 
Finally, in Section \ref{sect: Linear loci, linear measures and geometric orbifolds}, we discuss geometric structures on manifolds and orbifolds, introduce real linear orbifolds and equilinear orbifolds. We end by proving the results announced in the introduction. 

We feel that there is a nice parallel between the analysis of orbit closures and the analysis of invariant measures and our paper is structured to emphasize this. In particular we alternate between orbit closure discussions and measure discussions.
This approach may not appeal to every reader. A reader who wishes may read only the orbit closure discussion or may read the orbit closure discussion first and then proceed to the measure discussion.

\subsection{Acknowledgments}

We would like to thank Alex Eskin, Amir Mohammadi, Carlos Matheus, Simion Filip for helpful discussions and Carlos Ospina and Jiajun Shi for helpful comments on the manuscript.

\subsection{Disclosure on the use of generative AI tools}

Generative AI tools have been used to remove typographical errors in the text and improve clarity of the first paragraph of the Introduction. 

\section{Preliminaries}\label{sec: prelims}
In this section, we recall the basics of strata, paying special attention to the unnormalized case. There are many approaches to these definitions, we follow \cite[Sect. 2]{BSW}. A recent introduction to the field is contained in \cite{Athreya-Masur}. 

In this section we also discuss the connection between $\SL$-invariant measures on $\HH^1$ and homogeneous measures for the action of $\GL$ on $\HH$. Given an $\SL$-invariant and $d\ge 0$ we construct a $d$-homogeneous measure on $\HH$. We show that a $\GL$-homogeneous measure on $\HH$ gives rise to an $\SL$-invariant measures on $\HH^1$ via the cone construction. We also relate measures on $\HH$ to those on $\HH_m$. For completeness, we give a proof of an often quoted result about subgroups of the mapping class group.

\subsection{Strata and period coordinates}\label{subsec: strata}

Let $S$ be a connected, compact orientable surface of genus $g\geq 1$. Let $\Sigma = \{\xi_1, \ldots, \xi_n \}$ be a finite subset of $S$. Let  $a_1, \ldots, a_n$ be non-negative integers with $\sum a_j = 2g-2$. Let $\HH(a_1, \ldots, a_n)$ be the corresponding stratum of translation surfaces. We let $\HH_{\mathrm{m}}(a_1, \ldots, a_n)$ denote the stratum of marked translation surfaces $f:S\to M$ and we let $\pi : \HH_{\mathrm{m}}(a_1, \ldots, a_n) \to \HH(a_1, \ldots, a_n)$ be the forgetful mapping which drops the marking. 

We write $\Mod(S, \Sigma)$ for the group of isotopy classes of orientation-preserving homeomorphisms of $S$ stabilizing $\Sigma$ and preserving the type of the singularity modulo isotopies of $S$ that leave the set $\Sigma$ fixed. Thus an element of $\Mod(S, \Sigma)$ can switch two singular points of the same order. We also consider the pure mapping class group, $\PMod(S,\Sigma)$, which fixes the  points in $\Sigma$. The group $\Mod(S, \Sigma)$ acts properly discontinuously on the right on $\HHm$ by precomposition and the quotient is $\HH$. A simple example is given by the space of flat tori with a one marked point,  $\HH(0)$, which is $\GL/\mathrm{SL}(2,\Z)$ and $\HHm(0)$ is $\GL$. 
Note that $\HHm(0)$ is not simply connected so $\HHm(0)$ is not the {\em universal} cover of $\HH$. 
Note also that for surfaces of higher genus $\HH$ has finitely many components while $\HHm$ has infinitely many components (see \cite{Calderon-Salter}).
 
We will typically denote an element of $\HH$ by the letter $q$ when we want to consider it as a point of $\HH$, and by the letter $M$ or $M_q$ when we want to consider its underlying topological or geometrical properties as a space in its own right. Points in $\HHm$ will typically be denoted by boldface letters such as  $\q$. 

We recall the definition of the developing map which is written $\dev: \HH_{\mathrm{m}} \to H^1(S, \Sigma; \R^2)$. 
For an oriented path $\gamma$ in $M_q$ which is either closed or has endpoints at singularities, let $\hol(M_q, \gamma) \df \left( \int_\gamma dx_q, \int_\gamma dy_q \right)$, where $dx_q$ and $dy_q$ are the 1-forms on $M_q$ inherited from the forms $dx$ and $dy$ on the plane. 
Given $\q \in \HH_{\mathrm{m}}$ represented by $f: S \to M_q$, where $M_q$ is a translation surface, we define $\dev(\q) \df f^*(\hol(M_q, \cdot))$. The developing map is also known in the literature as the \textit{period map}. The space $\HH_{\mathrm{m}}$ has an atlas of charts given by an open cover $\{\mathcal{U}_\tau\}$, indexed by triangulations $\tau$ of $S$ with triangles whose vertices are in $\Sigma$. (Since we allow edges that do not have distinct vertices, the $\tau$ are not exactly triangulations. See \cite[Subsec. 2.4]{BSW}  for more details.)

The  maps
$$\phi_\tau \df \dev|_{\mathcal{U}_\tau}, \ \ \phi_\tau: \mathcal{U}_{\tau} \to H^1(S, \Sigma; \R^2)$$ are homeomorphisms onto their images. The charts $\phi_\tau$ give an atlas with linear overlap maps and endow $\mathcal{H}_{\mathrm{m}}$ with the structure of a manifold. (We will return to a discussion of this atlas in \S5.) This atlas of charts $\{(\mathcal{U}_\tau, \phi_\tau)\}$ is known as the {\em period coordinate atlas} on $\HHm$. 

Since the $\Mod(S,\Sigma)$ action on $\HH_{\mathrm{m}}$ is properly discontinuous and $\Mod(S,\Sigma)$ contains torsion elements, $\HH$ inherits the structure of an orbifold. In addition to the space $\HH$ it will be important to consider quotients of $\HHm$ by finite index subgroups of $\Mod(S,\Sigma)$. We will follow \cite{EM} and refer to the spaces $\HH_\Gamma$ as covers of $\HH$. For these spaces  to be literally covers we should interpret them as covering spaces in the sense of orbifolds (see Section \ref{sect: Linear loci, linear measures and geometric orbifolds} for more on this viewpoint). Let $\Gamma \subset \Mod(S,\Sigma)$ be such a group. We denote by $\HH_\Gamma$ the quotient of $\HHm$ by $\Gamma$, we denote by $\pi_\Gamma: \HHm \to \HH_\Gamma$ and by $p_\Gamma : \HH_\Gamma \to \HH$ the canonical projections. Note that we have $\pi = p_\Gamma \circ \pi_\Gamma$. 

An example of such a $\Gamma$ is given by $\PMod(S, \Sigma)$, the pure mapping class group, which is the finite-index normal subgroup of $\Mod(S, \Sigma)$ fixing $\Sigma$ pointwise. In this case we can identify $\HH_\Gamma$ with the space of translation surfaces with labeled singular points. There are different conventions in use in the literature regarding labeling of singular points of translation surfaces. All of our results hold in strata of translation surfaces with either labeled or unlabeled singularities. 

If $\Gamma$ is torsion free then the space $\HH_\Gamma$ inherits the structure of a manifold from $\HH$, in such a way that $\pi_\Gamma$ is a local diffeomorphism. 
If $\Gamma$ has torsion then the space $\HH_\Gamma$ is an orbifold. The discussion of period coordinates is a little more complicated in the setting of orbifolds. In this case each point in $\HH_\Gamma$ has a neighborhood, $N$, which is modeled on a finite quotient of an open set $\U\subset H^1(S,\Sigma,\R^2)$ by a linear action and it is convenient to use ``inverse charts" which map from $U$ to $N$ to define period coordinates. This is discussed in more detail in Section \ref{sect: Linear loci, linear measures and geometric orbifolds} but we will also give a quick summary here.

We denote by $\Gamma(\q) \subset \Gamma$ the stabilizer in $\Gamma$ of $\q \in \HHm$ and let $\U \subset \HHm$ be a $\Gamma(\q)$-invariant neighborhood of $\q$. We can assume that $\U$ is small enough so that $\dev|_\U$ is invertible and $\pi_\Gamma|_\U$ induces a homeomorphism $\U_{\Gamma(\q)} \to \pi_\Gamma(\U)$, where $\U_{\Gamma(\q)}$ denotes the quotient of $\U$ by $\Gamma(\q)$. Let $U = \dev(\U)$ and let $\varphi = \pi_\Gamma \circ (\dev|_\U)^{-1}$. The triple $(U,\Gamma(\q),\varphi)$ is a chart around $q$, that we call an {\em orbifold period chart}. 
Note that when $\Gamma$ is torsion free, then $\Gamma(\q)$ is trivial.   

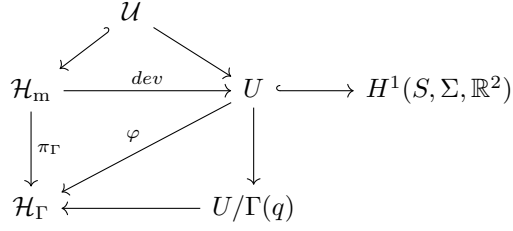
\begin{figure}
    \begin{tikzcd}[row sep=scriptsize, column sep=scriptsize]
        &\U \arrow[dr] \arrow[hook',dl]\\
        \HHm  \arrow[rr, "dev"] \arrow[dd, "\pi_\Gamma"] && U \arrow[dd] \arrow[ddll, "\varphi" '] \arrow[hook,r]& H^1(S,\Sigma,\R^2)\\
        \\
        \HH_\Gamma  && U/\Gamma(q) \arrow[ll]
    \end{tikzcd}
    \caption{Orbifold period charts for $\HH_\Gamma$}
    \label{fig: orbifold charts 1}
\end{figure}

The results of \cite{EM} deal with finite covers of $\HH$ which are manifolds. These can be constructed as spaces $\HH_\Gamma$ where $\Gamma$ is a torsion free subgroup of finite index in $\Mod(S,\Sigma)$.
Thus it is important for the results of \cite{EM} that manifold covers exist. This existence result is stated but not proved in \cite{EM}. For completeness we include a proof here.

\begin{thm}\label{thm: there exists a finite index normal subgroup}
    There exists a finite index, torsion free, normal subgroup $\Gamma\subset \Mod(S,\Sigma)$. 
\end{thm}

If $\Gamma$ is a subgroup of $\Mod(S,\Sigma)$ as in the statement of Theorem \ref{thm: there exists a finite index normal subgroup}, then $\HH_\Gamma$ is a finite manifold cover of the stratum $\HH$. If we take $\Gamma$ to be the intersection of this group with $\PMod(S,\Sigma)$ then $\HH_\Gamma$ is a finite manifold cover of the stratum with labeled singularities.

\begin{proof}
    We will give a proof in the case when $S$ is a surface of genus at least one. To begin with assume that $g\ge2$.
    Let $\Sigma_n \subset S$ be a subset of $S$  consisting of $n$ labeled points $\xi_1\dots \xi_n$. Let $\Sigma_0$ be the empty set. Let $\Mod(S,\Sigma_n)[3]$ denote the kernel of the homomorphism from $\Mod(S,\Sigma_n)$ to $Aut\big(H_1(S;\mathbb{Z}/3\mathbb{Z})\big)$ given by the homology representation of the mapping class group. This is a normal subgroup with finite index since $H_1(S;\mathbb{Z}/3\mathbb{Z})$ is finite.  Let $\Gamma_n$ denote the intersection of $\Mod(S,\Sigma_n)[3]$ with $\PMod(S,\Sigma_n)$. Since $\Gamma_n$ is the intersection of finite-index normal subgroups it is a finite index normal subgroup of $\Mod(S,\Sigma_n)$. 
    
    We show by induction on $\ell$ that $\Gamma_\ell$ is a torsion free subgroup of $\Mod(S,\Sigma_\ell)$. The case $\ell = 0$ is proved in \cite[Cor 1.8]{Ivanov}. Assume the result holds for a given $\ell$. Working in the pure mapping class group  the Birman exact sequence \cite[p. 98]{FM} gives: 
    $$
        1\rightarrow\pi_1(S-\Sigma_\ell,\xi_{\ell+1})\rightarrow \PMod(S,\Sigma_{\ell+1}) \rightarrow \PMod(S,\Sigma_\ell)\rightarrow 1,
    $$
    where the second map is the point pushing map and the third map is the forgetful map which ``forgets" the point $\xi_{\ell+1}$.
    If we restrict to the group $\Gamma_{\ell+1}$ we get an exact sequence 
    $$
        1\rightarrow\pi_1(S-\Sigma_\ell,\xi_{\ell+1})\rightarrow \Gamma_{\ell+1} \rightarrow \Gamma_{\ell}\rightarrow 1.
    $$
    If an element $\gamma$ of $\Gamma_{\ell+1}$ is torsion then by the induction hypothesis it maps to the identity in $\Gamma_\ell$ so it comes from an element of $\pi_1(S-\Sigma_\ell,\xi_{\ell+1})$. The fundamental group of an oriented surface is torsion free so $\gamma$ is the identity. In particular, $\Gamma_{\ell+1}$ is a finite-index, torsion free, normal subgroup of $\Mod(S,\Sigma)$. 

    When $g=1$ and $\ell=0$ the Birman exact sequence is not valid.  In this case the forgetting map is an isomorphism. Except for this point, the proof goes through as before.
\end{proof}

\subsection{The \texorpdfstring{$\GL$}{GL}-action and the \texorpdfstring{$\SL$}{SL}-action}\label{sec: from H to H1} 

The group $\GL$ acts on the left on the set of translation surfaces in $\HHm$ by modifying the planar charts of a translation surface. The group $\GL$ acts on $H^1(S, \Sigma; \R^2)$ via its action on the coefficients $\R^2$ and the  developing map $\dev:\HHm \to H^1(S,\Sigma;\R^2)$ is equivariant.  The sets $\mathcal{U}_\tau$ are invariant under this action. The following elements of $\SL$ will play an important role in dynamical questions: 

\begin{equation}\label{eq:usandgt}
    u_s = \left(\begin{matrix}  1 & s \\  0 & 1 \end{matrix} \right) \ \ \text{and} \ \ 
    g_t = \left(\begin{matrix}  e^t & 0 \\  0 & e^{-t} \end{matrix} \right).
\end{equation}
The action of the one-parameter subgroup $\{u_s \ | \ s \in \R\} \subset \SL$ is known as the {\em horocycle flow} and the one-parameter subgroup $\{g_t \ | \ t \in \R\}$ is the {\em geodesic flow}.

Consider the ``scalar multiplication" action of $\R_+^* $ on $\HH$ where the action of $t>0$ corresponds to the $\GL$ action given by the matrix
\begin{equation}
    \left(\begin{matrix}  t & 0 \\  0 & t \end{matrix} \right).
\end{equation}
Let $C$ denote the group of such matrices. The group $\GL$ is the product of $\SL$ and $C$ and, using scalar multiplication, we can write $\rho : \R^+_* \times \SL \to \GL$ for the natural isomorphism. 

The left $\GL$-action on $\HHm$ commutes with the right $\Mod(S, \Sigma)$-action on $\HHm$, thus $\GL$ also acts on $\HH$ in such a way that the map $\pi:\HHm\to \HH$ is $\GL$-equivariant. If $\Gamma \subset \Mod(S,\Sigma)$ is a finite index subgroup, we can also consider the induced action of $\GL$ on $\HH_\Gamma$. The canonical projections $\pi_\Gamma : \HH_{\mathrm{m}} \to \HH_\Gamma$ and $p_\Gamma : \HH_\Gamma \to \HH$ are $\GL$-equivariant. We denote by $\HH^1$ the set of area 1 surfaces in $\HH$, and we denote by $\HH_\Gamma^1$ the set of area 1 surfaces in $\HH_\Gamma$. The $\GL$ action on $\HH$ and $\HH_\Gamma$ restricts to an  $\SL$ action on $\HH^1$ and $\HH_\Gamma^1$. 

Many of the natural geometric objects that one wishes to consider on $\HH$ and on $\GL$ orbit closures in $\HH$ are homogeneous with respect to the scalar action. The natural measures that one wishes to consider in $\HH$ are typically not invariant under the action of $\GL$ but they do satisfy the following homogeneity property.

\begin{defin}\label{defin: homogeneous measure}
    Let $X$ be a topological space endowed with an action of $\GL$ and let $\delta >0$. A Radon measure $m$ on $X$ is said to be {\em homogeneous of degree $\delta$}, or simply {\em $\delta$-homogeneous}, if 
    $$
        \forall g \in \GL, \ g_\ast m = \mathrm{det}^{-\delta/2}(g) \cdot m. 
    $$ 
    The measure $m$ is said to be {\em ergodic} if every $\GL$-invariant Borel set of positive measure is conull. 
\end{defin}

For example a $\GL$-{\em invariant} measure is a homogeneous measure of degree zero. 

The notion of a homogeneous measure is inspired by the definition of a homogeneous function on a real vector space. It is also useful to consider the scalar multiplication action on measures on vector spaces. If $s_r$ denotes the scalar multiplication action of $\R_+^*$ on a vector space $V$ then we say a measure $m$ on $V$ is homogeneous of degree $\delta$ if 
$$
    (s_r)_* m=r^{-\delta} m.
$$
The factor of $2$ in Definition \ref{defin: homogeneous measure} is chosen to make these two notions compatible.

We denote by $\alpha : \HH \to \R_+^*$ the area function which is homogeneous of degree two. 
A computation shows that $g^*\alpha = \det(g) \cdot \alpha$ which is consistent with Definition \ref{defin: homogeneous measure}.

There is a basic connection between $\SL$ invariant measures on $\HH^1$ and homogeneous measures on $\HH$.  We describe a construction here which starts with an $\SL$ invariant measure on $\HH^1$ and produces a homogeneous measure on $\HH$. In Proposition \ref{prop : disintegration with respect to area function} we will give a construction which goes in the opposite direction.

We define the following map: 

\begin{align*}
    \Phi : \R_*^+ \times \HH^1 \to \HH, \ (t,q) \mapsto t\cdot q. 
\end{align*}
The map $\Phi$ is a bijection.  One sees that the map $\Phi$ intertwines the product action $\R^+_* \times \SL \curvearrowright \R_*^+ \times \HH^1$ with that of $\GL$ on $\HH$ via $\rho$. Since the action of $\R^+_*$ on itself by left multiplication is easy to describe, it is indeed enough to understand the action of $\SL$ on $\HH^1$. For any $\delta \in \R$, we let $\mu_\delta$ be the measure on $\R^+_*$ 
defined by the 1-form $t^{\delta-1} dt$. We have:

$$
    \forall f \in C_c(\R_*^+), \ \mu_\delta(f) = \int_{\R_*^+} f(t) t^{\delta-1} dt. 
$$
The measure $\mu_\delta$ satisfies $(s_r)_* \mu_\delta = r^{-\delta}\mu_\delta$ for all $r>0$. 

For any measure $\nu$ on $\HH^1$, we denote by $\mu_{\delta} * \nu$ the convolution of $\nu$ by $\mu_\delta$, \textit{i.e} the unique Radon measure on $\HH$ such that: 

\begin{align}\label{eq : disintegration with respect to area function}
    \forall f \in C_c(\HH), \ \mu_{\delta} * \nu(f) = \int_{\R^+_*} \left( \int_{\HH^1} f(t q) t^{\delta-1} d\nu(q) \right) dt. 
\end{align}

The measure $\mu_{\delta} * \nu$ is a homogeneous measure on $\HH$ of degree $\delta$.
In the terminology of Quint \cite{Quint} the measure $\mu_{\delta} * \nu$ would be denoted by $t^{\delta-1}dt\otimes \nu$.

\begin{prop}\label{prop : disintegration with respect to area function}
    Let $m$ be a homogeneous measure on $\HH$ of degree $\delta>0$. There is a finite $\SL$-invariant measure $\nu$ on $\HH^1$ such that $m=\mu_\delta * \nu$. The measure $\nu$ is uniquely determined by $m$ and $\nu$ is ergodic if and only if $m$ is. 
\end{prop}

\begin{rem}
    Note that the measure $\nu$ produced by the proposition is finite, even though the initial measure $m$ is not.
    It is the finiteness of $\nu$ which allows one to apply the tools of ergodic theory.
\end{rem}

\begin{proof} 
    We first prove that the pushforward measure $\alpha_\ast m$ of the measure $m$ by the area map $\alpha$ is a locally finite measure on $\R_+^\ast$. Let $K_1$ be a compact set in $\R_+^\ast$ and let $K_2 \subset \HH$ be a compact subset such that for any $q \in \alpha^{-1}(K_1)$ without horizontal saddle connections of length smaller than $1$, the following holds: 
    
    \begin{equation}\label{eq: nondivergence}
          \limsup \frac{1}{L} \big\lvert \{\ell \in [0,L] \ :\ u_\ell q \notin K_2 \} \big\lvert < \frac{1}{2}. 
    \end{equation}
    Such a $K_2$ exists by non-divergence for the horocycle flow, see \cite[Theorem 6.3]{minsky2002nondivergence}. The Birkhoff Ergodic Theorem for locally finite measures, see for example \cite[Theorem 1.14]{walterintroduction}, implies that there is a $U$-invariant function $f \in L^1(m)$ such that for $m$-a.e $q \in \HH$, we have
    
    \begin{equation}\label{eq : birkhoff}
         \frac{1}{L} \int_0^L 1_{K_2}(u_\ell q) \ d\ell  \underset{L \to \infty}{\to} f(q). 
    \end{equation}
    
    Let $\HH^{\geq t}$ be the set of surfaces in $\HH$ without horizontal saddle connections shorter than $e^{-t}$. Equations \eqref{eq: nondivergence} and \eqref{eq : birkhoff} imply that the function $f$ is $m$-a.s bounded below by $1/2$ on $\alpha^{-1}(K_1) \cap \HH^{\geq 0}$ and thus, for $f$ to be integrable, that 
    
    $$
        m(\alpha^{-1}(K_1) \cap \HH^{\geq 0}) \leq 2 \| f\|_{L^1(m)}.
    $$ 
    
    Now, using the fact that $g_t \HH^{\geq t} = \HH^{\geq 0}$, that $m$ is $g_t$-invariant and that the $\HH^{\geq t}$ is a nested sequence whose union is $\HH$, we obtain that $\alpha_{\ast}m(K_1) = m(\alpha^{-1}(K_1))$ is also bounded by $2 \| f\|_{L^1(m)}$. This shows that $\alpha_\ast m$ is locally finite. 
    
    Homogeneity of $m$ implies that $\alpha_\ast m$  is multiplied by $t^{-\delta}$ when rescaled by $t > 0$.  This uniquely determines the  measure $\mu_\delta$, up to a multiplicative constant. We can thus disintegrate $m$ on the fibers of the map $\alpha$ with respect to $\mu_\delta$ and we get 
    $$
        m = \int_{\R_+^*} m_t \ d\mu_\delta(t).  
    $$
    Uniqueness of fiberwise measures and homogeneity of $m$ under scalar matrices imply that for any $t>0$ and for a conull set of $s>0$, we have $t_\ast m_s = m_{st}$. A Fubini-type argument shows that there is $s_0 >0$ such that the previous equality holds for a conull set of $t>0$ and we let $\nu = (1/s_0)_\ast m_{s_0}$. We compute: 
    
    \begin{align*}
        \mu_\delta * \nu(f) 
        &= \int_{\R_+^*} \int_{\HH^1} f(q) t^{\delta-1} d t_\ast \nu(q)dt \\
        &= \int_{\R_+^*} \int_{\HH^1} f(q) t^{\delta-1} d m_t(q)dt
        = m(f). 
    \end{align*}
    It is clear from the above considerations that the measure $\nu$ is finite, uniquely determined by equation \eqref{eq : disintegration with respect to area function} and in particular $\SL$-invariant. The last statement follows from equation \eqref{eq : disintegration with respect to area function} and the observation that a Borel set $A\subset\HH$ is $\GL$-invariant if and only if it has the form $A = \R_*^+A_1$ where $A_1$ is an $\SL$-invariant subset of $\HH^1$.   
\end{proof}

\begin{rem}
    Notice that we only used the fact that $m$ is homogeneous under the action of the upper triangular subgroup of $\GL$ for the existence and uniqueness of the measure $\nu$ in the proof of Proposition \ref{prop : disintegration with respect to area function}.
\end{rem}

In Proposition \ref{prop : coning is the same} we will see that the measure $\nu$ appearing in the previous proposition can also be described as the cone measure as in \cite{Athreya-Masur}. 
We first introduce some terminology. If $B$ is any Borel set in $\HH^1$, we let $C(B) := \{tq \ : \ q \in B, \ 0<t \leq 1\}$ be the unit cone over $B$ in $\HH$. If $m$ is any  measure on $\HH$, we define a measure $m_1$ on $\HH^1$ by the formula $m_1(\cdot) = m(C(\cdot))$. This measure is called the {\em cone measure} associated to $m$. Note that the cone measure could \textit{a priori} be everywhere infinite. However, we have the following : 

\begin{prop}\label{prop : coning is the same}
    Let $m$ be a homogeneous measure on $\HH$. The corresponding cone measure $m_1$ is finite and coincides, up to a multiplicative constant, with the measure $\nu$ appearing in Proposition \ref{prop : disintegration with respect to area function}. 
\end{prop}

\begin{proof}
    We will prove both statements simultaneously. Let $\delta$ be the degree of $m$, let $\nu$ be as in Proposition \ref{prop : disintegration with respect to area function} and let $B \subset \HH^1$ be a Borel set.  We compute

    \begin{align*}
        m_1(B) = m(C(B)) &= \int_{\R_*^+} \int_{\HH^1} 1_{C(B)}(tq) t^{\delta-1} d\nu(q) dt \\  &= \int_{\HH^1} 1_{C(B)}  \Bigl( \int_0^1 t^{\delta-1} \ dt \Bigr)d\nu(q) = \frac{1}{\delta}\nu(B).
    \end{align*}
\end{proof}

\subsection{Pre-image measures}\label{quotientmaps}
In \cite{SSWY} we describe the pre-image construction which will play a role here.
Let $X$ be a manifold and let $\Lambda$ be a group acting properly discontinuously on the right. We write $X_\Lambda$ for the quotient space and $p : X \to X_\Lambda$ for the \textit{quotient map}. We denote by $\Lambda(x) \subset \Lambda$ the stabilizer of $x \in X$. To a Radon measure $m$ on $X_\Lambda$ we associate a $\Lambda$-invariant Radon measure $\widetilde m$ on $X$. This measure is characterized by the fact that it satisfies: 
$$
\forall f \in C_c(X), \ \int_X f \ d\widetilde m = \int_{X_\Lambda} \sum_{x \in p^{-1}(q)} f(x)\big|\Lambda(x)\big| dm(q)
$$
in light of the Riesz Representation Theorem.
We call the measure $\widetilde m$ the pre-image of $m$ by the quotient map $p$, and we denote it by $p^*m$. 

Every $\Lambda$-invariant Radon measure on $X$ is the pre-image measure for a unique measure $m$. See \cite[Appendix B]{SSWY} for more details. Examples to which this result can be applied include the case where $X = \HHm$ and $\Lambda=\Gamma$ is a subgroup of $\Mod(S,\Sigma)$. When $\Gamma \subset \Mod(S,\Sigma)$ is a normal subgroup, we can consider the case where $X = \HH_\Gamma$ with the action of $\Gamma_0 = \Gamma \backslash \Mod(S,\Sigma)$ (note that $\HH_\Gamma / \Gamma_0 \simeq \HH$). We will also consider the case where $X = U$ and $\Lambda = \Delta$ where $(U,\Delta,\varphi)$ is an orbifold period chart for $\HH_\Gamma$.  

In the case where $\Lambda$ is finite, the image of a measure is closely related to its push-forward: 

\begin{prop}\label{prop: image of a measure in the finite case}
    If $\Lambda$ is finite, then $m = \frac{1}{|\Lambda|} \cdot p_\ast \widetilde{m}$.  
\end{prop}

\begin{proof}
    Let $f \in C_c(X)$. Notice that if $q = p(x_0)$, then $\sum_{x \in p^{-1}(q)} f(x)\big|\Lambda(x)\big| = \sum_{\lambda \in \Lambda} f(x_0\cdot \lambda)$ and thus
        \begin{align*}
            \frac{1}{|\Lambda|}\int_{X_\Lambda} \sum_{x \in p^{-1}(q)} f(x)\big|\Lambda(x)\big| dp_*\widetilde{m}(q) &= \frac{1}{|\Lambda|}\int_{X} \sum_{\lambda \in \Lambda} f(x_0\cdot \lambda) d\widetilde m(x) \\
            &= \int_X f \ d\widetilde m. 
        \end{align*}
    The last equality uses the $\Lambda$ invariance of $\widetilde{m}$. This ``disintegration-like'' property for $\widetilde m$ uniquely determines $m$, see \cite[Prop. B3]{SSWY} for a precise statement. 
\end{proof}

\section{Black boxing some results of Eskin, Mirzakhani and Mohammadi}\label{sec: recasting}

We start by expanding on the discussion in the Introduction. While passing to a finite cover does not change the local property of orbit closures, generic injectivity of the parameterizing map could be lost in the process. 
To see this one could use the results of \cite{EM} and \cite{EMM} for finite covers to construct a parametrization of an orbit closure in a stratum. 
Say $\M\subset\HH$ is an $\SL$ orbit closure. We could pick a ``cover" $\widetilde\HH$ with no orbifold points, pick an orbit closure $\widetilde \M \subset \widetilde\HH$ mapping to $\M\subset\HH$ by the ``covering" map $p:\widetilde\HH\to \HH$ and choose a parametrization $f':\N\to\widetilde\HH$ of $\widetilde \M$ using the results of \cite{EMM}. The composition $p \circ f'$ is then a parametrization of $\M$. It is indeed proper and ``locally linear'' but it need not be generically injective. 

Our work relies on the work of \cite{EM} and \cite{EMM} for the description of the local properties of orbit closures and invariant measures in finite manifold covers to prove global properties of orbit closures and invariant measures in strata. We start by giving names to the properties we will use.

We will use the fact that $\GL$-orbit closures are locally given as unions of finite collections of linear subspaces in period coordinates. We formalize this property with the following definition. Recall that period charts are defined in subsection \ref{subsec: strata}. 

\begin{defin}\label{def: locally linear}
    A {\em locally linear set} in $\HH_\Gamma$ is a closed set $\M \subset \HH_\Gamma$ such that for any period chart $(U,\Delta,\varphi)$, there are finitely many linear subspaces $(V_i)_{i \in I} \subset H^1(S,\Sigma;\C)$ defined over $\R$ such that: 
    \begin{equation}\label{eq: finitely many linear pieces}
       \varphi^{-1}(\M) = U \cap \bigcup_{i\in I} V_i.  
    \end{equation}
\end{defin}

The fact that we have $\varphi^{-1}$ is a result of the fact that we are using orbifold charts which is appropriate when $\Gamma$ has torsion and $\HH_\Gamma$ is an orbifold. (See the discussion before subsection \ref{thm: there exists a finite index normal subgroup} for more on orbifold charts.) If $\Gamma$ is torsion free then $\HH_\Gamma$ is a manifold. In this case we have charts $\varphi$ which map in the usual direction $\varphi:U_\tau\to U\subset H^1(S,\Sigma;\R^2)$ and the condition would be that 
\begin{equation}
       \varphi(\M\cap U_\tau) = U \cap \bigcup_{i\in I} V_i.  
\end{equation}

We will also use the fact that ergodic $\SL$-invariant measures on $\HH_\Gamma^{1}$ come from measures on $\HH_\Gamma$ that are locally given as a sums of finite collections of Lebesgue measures on linear subspaces in period coordinates. We formalize this property with our notion {\em locally linear measure} (compare Definition 1.1 in \cite{EM}).

\begin{defin}\label{def: locally linear measure}
    A {\em locally linear measure} on $\HH_\Gamma$ is a measure $m$ such that for any period chart $(U,\Delta,\varphi)$, there exists a finite collection of linear subspaces $(V_i)_{i\in I} \subset H^1(S,\Sigma;\C)$ defined over $\R$ such that: 
    \begin{enumerate}
        \item the pre-image measure $\varphi^*m$ is supported on $U \cap \bigcup_{i\in I} V_i$,
        \item for all $i\in I$, the measure $\varphi^* m|_{V_i}$  is the restriction of a Lebesgue measure of $V_i$.
    \end{enumerate} 
\end{defin}

In Definition \ref{def: locally linear measure} above, we use the pre-image construction for measures which is described in subsection \ref{quotientmaps}. There is a slight abuse of notation in the first item: by $\varphi^*m$, we mean the pre-image by the quotient map $U \to U/\Delta$ of $(\varphi_\Delta^{-1})_*m$, where $\varphi_\Delta : U_\Delta \to \varphi(U)$ is the homeomorphism induced by $\varphi$. This is the measure on $U$ which is uniquely defined by the property that for any compactly supported continuous function: 
$$
    \int_U f \ d\varphi^*m = \int_{\varphi(U)} \sum_{x \in \varphi^{-1}(p)}|\Delta(x)| f(x) \ dm(p). 
$$
Notice in particular that a locally linear measure is supported on a locally linear set. 

As discussed in the introduction, \cite{EM} and \cite{EMM} work in the setting of spaces $\HH_\Gamma$ where $\Gamma \subset \Mod(S,\Sigma)$ is torsion free (recall Thm. \ref{thm: there exists a finite index normal subgroup}). The results that we use are contained in the following statements.

\begin{thm}[\cite{EMM}]\label{thm: emm}
    Let $\Gamma \subset \Mod(S,\Sigma)$ be a finite index torsion free subgroup. The $\GL$-orbit closure of any $q \in \HH_\Gamma$ is a connected locally linear set.
\end{thm}

Note that the local finiteness property is not mentioned explicitly in the theorems stated in \cite{EM} and \cite{EMM}. We could deduce local finiteness from the properness of the map $f$ which is mentioned in Definition 1.1 of \cite{EM}. Local finiteness also follows from the fact that orbit closures are varieties \cite{varieties} but local finiteness is used as a {\em hypothesis} for this result.

In fact \cite{EM} does not deduce local finiteness from properness but proves local finiteness directly. 
The proof uses the existence of a probability measure $\nu$ supported on $\M^1$.
The idea of the proof is to show that in a neighborhood of a point $q\in\M^1$ each local sheet makes a definite contribution to the $\nu$ measure of $\M^1$. Since the total mass of $\nu$ is finite the number of local sheets must be finite. The measure $\nu$ on a sheet is originally described in terms of the restriction of a symplectic form to the tangent space of a sheet but this description does not give a positive lower bound for the mass of a sheet. Lemma 14.1 of \cite{EM} shows that this measure can also be described in terms of the restriction of the Hodge norm and this description does give a positive lower bound for the mass. The statement of local finiteness is at the bottom of page 273 of \cite{EM}.

\begin{thm}[\cite{EM}]\label{thm: em}
    Let $\Gamma \subset \Mod(S,\Sigma)$ be a finite index torsion free subgroup. Let $\nu$ be an ergodic  $\SL$-invariant probability measure on $\HH_\Gamma^{1}$. Then there is $\delta > 0$ such that $\lambda=\mu_\delta * \nu$ is a locally linear measure on $\HH_\Gamma$. 
\end{thm}

In light of Proposition \ref{prop : coning is the same}, the measure $\nu$ in Theorem \ref{thm: em} above can equivalently be described as the cone measure associated to a locally linear measure. 

Prior to Definition 1.1 in \cite{EM} the measure $\lambda$ is described as $d\lambda=c d\nu d\alpha$ where $\alpha$, the area function, is a homogeneous function of degree 2.
We read this as suggesting that $\lambda$ is a homogeneous measure of degree 2. This would not be correct.

Another approach to this result is to observe that towards the end of the proof of Lemma 16.3 in \cite{EM} 
there is an assertion that the support of the invariant measure is contained in finitely many ``affine pieces". This shows that $\lambda$ is horospherical. Combining this result with \cite[Theorem 1.2]{SSWY} shows that $\lambda$ is built from the cone construction. 

We will use these results in the proofs of Theorem \ref{thm: orbit is linear} and Theorem \ref{thm: invariant measure is cone measure associated to linear measure}. Theorem \ref{thm: orbit is linear} will be used in the proofs of Theorem \ref{thm: linear structure on orbit closures}. Theorem \ref{thm: invariant measure is cone measure associated to linear measure} will be used in the proof of Theorem \ref{thm: linear structure on invariant measures}.

\section{Linear loci and linear measures}\label{sec: linear loci and linear measures}
The paper \cite{SSWY} used a version of the classic argument from Margulis' thesis which requires the construction of product boxes. These boxes are initially constructed in $\HHm$ before being projected to $\HH$. For this purpose, the most convenient description of a $\GL$-orbit closure was a description of its pre-image by the map $\pi : \HH_{\mathrm{m}} \to \HH$. This description is contained in the definition of a linear locus (Definition \ref{defin: linear locus}). In Proposition \ref{prop: locally linear locus is linear} we prove that the local property of being a locally linear set implies the global property of being a linear locus. 

We introduce the notion of {\em linear measure} to describe a certain global geometric structure for measures in $\HH$.  We show that locally linear measures in $\HH$ are linear measures (Proposition \ref{prop: locally linear measure is linear}). 
We also introduce the notion of {\em irreducibility} for linear loci and linear measures. We show that a linear locus is a union of irreducible linear loci.  We deduce from the results of Section \ref{sec: recasting} that a $\GL$-orbit closure is an irreducible linear locus (Theorem \ref{thm: orbit is linear}). We show that an ergodic $\SL$-invariant measure is precisely the cone measure associated to an irreducible linear measure (Theorem \ref{thm: invariant measure is cone measure associated to linear measure}).

\subsection{Linear loci in \texorpdfstring{$\HH$}{HH}}

\begin{defin}\label{defin: real structure} A {\em real structure} on a complex vector space $V$ is given by the choice of a real linear subspace $W\subset V$ so that $W\cap iW=\{0\}$ and $V=W\oplus iW$.
\end{defin}

The space $V=H^1(S,\Sigma;\C)$ is an example of a vector space with a real structure where we take $W$ to be $H^1(S,\Sigma;\R)$.

\begin{defin}\label{defin: linear submanifold}
    A {\em linear submanifold} in $\HHm$ is a connected component of $\dev^{-1}(V)$ where $V$  is a linear subspace of $H^1(S,\Sigma;\R^2)$ with a real structure.
\end{defin}

Notice that a linear submanifold is indeed a manifold, this is a consequence of the fact that $\dev$ is a locally affine. If $L$ is a linear submanifold, we will usually denote by $V_L$ the linear subspace appearing in Definition \ref{defin: linear submanifold}, and we shall say that $L$ is modeled on $V_L$. 

\begin{defin}\label{defin: linear locus}
    A {\em linear locus} is a subset  $\M \subset \HH$ such that $\pi^{-1}(\M)$ is the union of a locally finite collection of linear submanifolds.  
\end{defin}

Recall that a collection $(L_i)_{i\in I}$ of subsets in $\HHm$ is said to be locally finite if around any point there is an open subset $\mathcal{U}$ which is intersected by at most finitely many of the $L_i$. Equivalently (because $\HHm$ is locally compact), this means that every compact set intersects at most finitely many of the $L_i$. In \cite[Definition 2.4]{SSWY} we used the term ``invariant subvariety'' for an object closely related to a linear locus. Linear loci are slightly more general than invariant subvarieties as no condition is imposed on the dimensions of the collection of spaces $V_L$.
In subsection \ref{subsec: decomposing linear loci into irreducible components} we will introduce {\em irreducible} linear loci and show they have a well-defined dimension. 

\begin{prop}\label{prop: invariant by the action of GL}
    A linear locus is closed and $\GL$-invariant.
\end{prop}

\begin{proof}
    Let $\M$ be a linear locus. We start by noting that a linear submanifold $L \subset \HH_{\mathrm{m}}$ is closed because it is a component of $\dev^{-1}(V_L)$ which is closed. The set $\pi^{-1}(\M)$ is closed because it is a locally finite union of closed sets, and this implies that $\M$ is closed.
    
    In order to show that $\M$ is $\GL$-invariant, it is enough to prove that $\pi^{-1}(\M)$ is $\GL$-invariant, since $\pi$ is $\GL$-equivariant.  Let $L$ be a linear submanifold contained in $\pi^{-1}(\M)$. By definition, $V_L$ is defined over $\R$ and by \cite[Lemma 2.1]{SSWY} it is invariant under the action of $\GL$ on $H^1(S,\Sigma;\R^2)$. Since $\dev$ is $\GL$-equivariant the action of $\GL$ on $\HHm$ preserves $\dev^{-1}(V_L)$. Since $\GL$ is connected, the action of $\GL$ on $\HHm$ preserves $L$. Since $\pi^{-1}(\M)$ is a union of linear submanifolds it follows that it is invariant under $\GL$.  
\end{proof}

\subsection{Decomposing linear loci into irreducible components}\label{subsec: decomposing linear loci into irreducible components}
We introduce the notion of {\em irreducibility} for linear loci and show that any linear locus can be canonically decomposed into a finite union of irreducible linear loci which we call {\em linear components}. 

\begin{defin}\label{def: linear locus}
    A linear locus is {\em irreducible} if it cannot be written as a union of two proper distinct linear subloci.  
\end{defin}

\begin{defin}
    An irreducible component of a linear locus, $\M$, is a maximal irreducible linear sublocus contained in $\M$. 
\end{defin}

\begin{lem}\label{lem: a lemma about inclusion of linear submanifolds} Let $L$ be a linear submanifold of $\HHm$.
    Let $(L_i)_{i \in I}$ be a countable collection of linear submanifolds in $\HHm$. If $\U$ is an open subset in $\HH_{\mathrm{m}}$ such that $L \cap \U \subset \bigcup_{i\in I} L_i$, then there is $i_0 \in I$ such that $L \subset L_{i_0}$. 
\end{lem}

\begin{proof}
    We have the equality $L \cap \U= \bigcup_{i \in I} L_i \cap L \cap \U$. It follows from the Baire category theorem that there is an $i_0 \in I$ such that $L_{i_0} \cap L \cap \U$ has non empty interior in $L \cap \U$. By the definition of a linear submanifold, this is possible only when $V_L \subset V_{L_{i_0}}$. This implies that $L \subset L_{i_0}$. 
\end{proof}

\begin{prop}\label{prop: max lin manifold}
    Let $\M$ be a linear locus. If $L$ is a maximal linear manifold in $\pi^{-1}(\M)$ then its image by $\pi$ is an irreducible component of $\M$. 
\end{prop}

\begin{proof}
    Denote by $\M_L$ the image of $L$ under the map $\pi$. We first prove that $\M_L$ is a linear locus. Let $(L_i)_{i \in I}$ be a locally finite collection of linear submanifolds whose union is $\pi^{-1}(\M)$. This is a countable collection because $\HHm$ is $\sigma$-finite. The pre-image of $\M_L$ by the map $\pi$ is the union of the translates of $L$ by elements of $\Mod(S,\Sigma)$. It follows from Lemma \ref{lem: a lemma about inclusion of linear submanifolds} in the case $U=\HH_{\mathrm{m}}$ that for every $\gamma \in \Mod(S,\Sigma)$, there is an index $i_\gamma \in I$ such that $L \cdot \gamma \subset L_{i_\gamma}$. Maximality of $L$ implies $L \cdot \gamma=L_{i_\gamma}$. Let $J = \{i_\gamma \ | \ \gamma \in \Mod(S,\Sigma)\} \subset I$. We just proved that $\pi^{-1}(\M_L) = \bigcup_{j \in J} L_{j}$. This proves that $\M_L$ is indeed a linear locus. 
    
    We must now show that $\M_L$ is irreducible and in order to do so,  it is enough to show that a linear locus $\N \subset \M_L$ with non-empty interior in $\M_L$ coincides with $\M_L$. Let $\N$ be such a linear locus and let $(L_j)_{j \in J}$ be a locally finite collection of linear submanifolds whose union is $\pi^{-1}(\N)$. Let $V$ be an open set in $\HH$ intersecting $\M_L$ such that $V \cap \M_L \subset \N$. This implies that $\pi^{-1}(\M_L) \cap \pi^{-1}(V) \subset \pi^{-1}(\N)$. Letting $\U = \pi^{-1}(V)$ in Lemma \ref{lem: a lemma about inclusion of linear submanifolds}, we obtain $\gamma \in \Mod(S,\Sigma)$ and $i \in I$ such that $L \cdot \gamma \subset L_i$. This implies $\M_L \subset \N$ and thus $\M_L = \N$. 
    
    It remains to show that $\M_L$ is maximal. Let $\M_0 \subset \M$ be an irreducible linear locus that contains $\M_L$. Let $(L_k)_{k \in K}$ be a locally finite collection of linear submanifolds whose union is $\pi^{-1}(\M_0)$. Letting $\U = \HH_{\mathrm{m}}$ in Lemma \ref{lem: a lemma about inclusion of linear submanifolds} shows that there is a $k_0 \in K$ such that $L \subset L_{k_0}$. Maximality of $L$ implies $L = L_{k_0}$. Let $\M_L^{\perp}$ be the image of the union of the $L_k$ which do not coincide with a translate of $L$ by the mapping class group. This is a linear locus, and we get the following decomposition: 
    $$
        \M_0 = \M_L \cup \M_L^{\perp}.  
    $$
    By irreducibility, we obtain that $\M_L^{\perp} = \emptyset$ and $\M_0 = \M_L$ as desired. 
\end{proof}

\begin{prop}\label{prop:decomposin linear loci}
    A linear locus $\M$ is the union of its irreducible components and there are finitely many such components. 
\end{prop}

\begin{proof}
    We first prove that $\M$ is covered by finitely many of its irreducible components. Let $K \subset \HH^1$ be a compact set such that for any $q \in \HH^1$ without horizontal saddle connections of length smaller than $1$, the following holds: 

    \begin{equation}\label{eq : a use of nondivergence}
        \limsup \frac{1}{L} \big\lvert \{\ell \in [0,L] \ :\ u_\ell q \notin K \} \big\lvert < \frac{1}{2}. 
    \end{equation}
    Such a $K$ exists by non-divergence for the horocycle flow, see \cite[Theorem 6.3]{minsky2002nondivergence}. Let $\tilde K \subset \HHm$ be a compact set such that $\pi(\tilde{K}) = K$. It follows from the definition of linear loci that there are finitely many linear submanifolds $L_1,\dots,L_\ell$ whose union contains $\pi^{-1}(\M) \cap \tilde K$. We can assume that the $L_i$ are maximal. Let $\M_i\subset \M$ denote the image of $L_i$. By Proposition \ref{prop: max lin manifold}, this is an irreducible linear sublocus. We have $K \cap \M \subset \bigcup \M_i$. Let $q \in \M$, let $\tau \in \R^+_*$ be such that $\tau \cdot q \in \HH^1$ and denote by $q_1$ the surface $\tau \cdot q$. The surface $g_t \cdot q_1$ does not have horizontal saddle connection of length shorter that $1$ if $t$ is large enough. By \eqref{eq : a use of nondivergence}, there is $s>0$ such that $u_s g_t \cdot q_1 \in K$. By Proposition \ref{prop: invariant by the action of GL}, $\M$ is $\GL$-invariant. This means that $u_s g_t \cdot q_1 \in K \cap \M$ and using again the fact that linear loci are invariant, we get that $q \in \bigcup \M_i$.  

    Now, let $\N \subset \M$ be an irreducible component of $\M$. We have the equality $\N = \bigcup \N \cap \M_i$. An intersection of linear loci is again a linear locus so by irreducibility, there is a index $i$ such that $\N = \N \cap \M_i$ and thus $\N \subset \M_i$. By maximality of $\N$, this shows that $\N = \M_i$.  
\end{proof}

Let $\M$ be an irreducible linear locus and let $L$ be a maximal linear submanifold in $\pi^{-1}(\M)$. It follows from Proposition \ref{prop: max lin manifold} and Proposition \ref{prop:decomposin linear loci} that $\pi(L)=\M$. We call such a linear submanifold a \textit{lift} of $\M$. It is easy to verify that any two lifts of $\M$ differ by the action of an element of $\Mod(S,\Sigma)$ and thus have same dimension.  

\begin{defin}
    The {\em (complex) dimension} of an irreducible linear locus $\M$ is the complex dimension of any of its lifts. We denote it by $\dim_\C(\M)$. 
\end{defin}

\subsection{Locally linear sets and linear loci in \texorpdfstring{$\HH$}{HH}}

In Section \ref{sec: recasting}, we have defined locally linear sets in $\HH_\Gamma$ for arbitrary $\Gamma$. In this section we focus on the case where $\HH_\Gamma$ is the stratum $\HH$
and we write $\HH$ for $\HH_\Gamma$, that is the case where $\Gamma$ is $\Mod(S,\Sigma)$.
In fact our proof also applies to $\HH_\Gamma$ where $\Gamma$ has finite index in $\Mod(S,\Sigma)$. In particular we can also take $\HH$ to be the unmarked stratum where $\Gamma=\PMod(S,\Sigma)$.

\begin{prop}\label{prop: locally linear locus is linear}
    A locally linear set in $\HH$ is a linear locus.
\end{prop}

We first state two preliminary lemmas that will be key to the proof of Proposition \ref{prop: locally linear locus is linear}. 

\begin{lem}\label{lem: locally contained implies contained}
    Let $\M$ be a locally linear set in $\HH$ and let $L$ be a linear submanifold in $\HHm$. If there exists an open set $\W \subset \HHm$ such that $\pi(\W \cap L) \subset \M$, then $\pi(L) \subset \M$. 
\end{lem}

\begin{proof}
    Let $\M$ and $L$ be as in the statement of the Lemma. Let $L_0$ be the set of points in $L$ that are contained in some neighborhood $\U \subset \HHm$ such that $\pi(\U \cap L) \subset \M$. The set $L_0$ is clearly open in $L$ and it is not empty by assumption, we are going to show it is also closed. Let $(\q_n )_n$ be a sequence in $L_0$ converging to $\q_\infty \in L$. Since $\M$ is closed, we know that $\pi(\q_\infty)$ is contained in $\M$. Let $\U \subset \HHm$ be a neighborhood of $\q_\infty$ and $\Delta \subset \Mod(S,\Sigma)$ be a finite subgroup such that $\pi$ induces a homeomorphism $\U_\Delta \to \pi(\U)$ and $\dev|_{\U}$ is invertible. Let $U=\dev(\U)$ and let $\varphi = \pi \circ (\dev|
    _{\U})^{-1}$. The triple $(U,\Delta,\varphi)$ is a period chart for a neighborhood $\pi(\q_\infty)$, as depicted below:
    
    \begin{figure}[h!]
        \begin{tikzcd}[row sep=scriptsize, column sep=scriptsize]
            &\U \arrow[dr] \arrow[hook',dl]\\
            \HHm  \arrow[rr, "dev"] \arrow[dd, "\pi"] && U \arrow[dd] \arrow[ddll, "\varphi" '] \arrow[hook,r]& H^1(S,\Sigma,\R^2)\\
            \\
            \HH  && U/\Delta \arrow[ll]
        \end{tikzcd}
        \label{fig: orbifold charts 2}
    \end{figure}

    By the definition of locally linear sets, there is a finite collection of linear subspaces $V_1,\dots,V_\ell \subset H^1(S,\Sigma;\R^2)$ defined over $\R$ such that 
    
    $$\varphi^{-1}(\M) = U \cap \bigcup_{i=1}^\ell  V_i.$$ 
    For every $n>0$, there is a non-empty open set $\U_n \subset \HHm$ around $\q_n $ such that the set $\pi(\U_n \cap L)$ is contained in $\M$. For $n$ large enough, we have $\q_n  \in \U$ and thus $\U \cap \U_n$ is a nonempty open subset. The set $\dev(\U \cap \U_n \cap L)$ is an open subset of $V_L$ which is contained in $\bigcup_{i=1}^\ell V_i$. Since a linear subspace is spanned by any of it open subset, we deduce that the entire $V_L$ is contained in that finite union. We conclude that the entire vector space $V_L$ is contained in the union of the $V_i$ and thus there is an index $i \in \{1,\dots,\ell\}$ such that $V_L \subset V_i$ as $V_L$ would be a finite union of nowhere dense subsets otherwise. This shows that $\pi(\U \cap \U_n \cap L) \subset \M$ and thus $\q_\infty \in L_0$, as desired. The set $L_0$ therefore coincides with $L$, as $L$ is connected. This implies in particular that $\pi(L)$ is contained in $\M$.  
\end{proof}

\begin{lem}\label{lem: the components are locally finite}
    For any linear subspace $V \subset H^1(S,\Sigma;\R^2)$, the collection of connected components of $\dev^{-1}(V)$ is locally finite.   
\end{lem}

\begin{proof}
    Let $\U$ be an open subset of $\HHm$ and $U$ a convex open subset of $H^1(S,\Sigma;\C)$ such that $\dev$ induces a homeomorphism from $\U$ to $U$. We claim that $\U$  meets at most one connected component of $\dev^{-1}(V)$. Indeed, if $\boldsymbol{q}_1$ and $\boldsymbol{q}_2$ are in $\dev^{-1}(V) \cap \U$, then they are connected by the path $t \mapsto \dev|_\U^{-1}\big( \dev(q_1) + t(\dev(q_2)-\dev(q_1)\big)$ which is contained in $\dev^{-1}(V)$. 
    To conclude, it is enough to notice that $\HHm$ can be covered by open sets $\U$ as above. 
\end{proof}

\begin{proof}[Proof of Proposition \ref{prop: locally linear locus is linear}] 
    Let $\M \subset \HH$ be a locally linear set. Let $\mathfrak{L}$ be the collection of linear submanifolds contained in $\pi^{-1}(\M)$. This collection is nonempty and covers $\pi^{-1}(\M)$, this is a consequence of the definition of a locally linear set combined with Lemma \ref{lem: locally contained implies contained}. The collection $\mathfrak{L}$ is partially ordered by inclusion. We denote by $\mathfrak{L}_0$ be the sub-collection of linear submanifolds that are maximal with respect to the partial order. (Note that the length of a chain in $\mathfrak{L}$ is bounded by $\dim_\C H^1(S,\Sigma;\R^2)$ and thus existence of maximal elements does not require Zorn's lemma.) We are going to prove that the collection $\mathfrak{L}_0$ is locally finite. Let $\U \subset \HHm$ be a relatively compact open neighborhood and $\Delta \subset \Mod(S,\Sigma)$ be a finite group such that $\pi$ induces a homeomorphism $\U_\Delta \to \pi(\U)$ and $\dev_{|\U}$ is invertible. Let $U = \dev(\U)$ and let $\varphi=\pi \circ (\dev|_\U)^{-1}$. The triple $(U,\Delta,\varphi)$ is a period chart. By definition, there are linear subspaces $V_1,\dots,V_\ell \subset H^1(S,\Sigma;\R^2)$ such that

    $$
        \varphi^{-1}(\M) = U \cap \bigcup_{i=1}^\ell V_i.
    $$
    
    Let $L$ be a linear submanifold in $\mathfrak{L}_0$ that intersects $\U$. In particular, the set $U \cap V_L$ is contained in $\bigcup_{i=1}^\ell V_i$. We deduce as in the proof of Lemma \ref{lem: locally contained implies contained} that there is $i \in \{1,\dots,\ell\}$ such that $V_L \subset V_i$. By maximality, we deduce that $L$ is a connected component of $\dev^{-1}(V_i)$. Since $\U$ is relatively compact, there are at most finitely many such components that intersect $\U$ by Lemma \ref{lem: the components are locally finite}. This proves that $\U$ meets finitely many linear submanifolds in $\mathfrak{L}_0$. This concludes the proof of the local finiteness of $\mathfrak{L}_0$ as $\U$ was arbitrary.  
\end{proof}

\subsection{Linear Measures}\label{sect: linear measures}

\begin{defin}
    A {\em linear measure on $\HHm$} is a measure supported on a linear submanifold $L$ that locally corresponds to a Lebesgue measure on $V_L$ via $\dev$. More precisely, there is a Lebesgue measure $\lambda$ on $V_L$ such that for any $\U \subset \HHm$ on restriction to which $\dev$ is invertible, we have $(\dev|_\U)_\ast m|_{\U} = \lambda|_{\dev(\U)}$.   
\end{defin}

\begin{defin}\label{def: linear measure}
    A measure $m$ on $\HH$ is called linear if there is a locally finite collection of linear submanifolds $(L_i)_{i \in I}$ in $\HH_{\mathrm{m}}$ whose union is the support of $\pi^* m$ and such that $\pi^*m|_{L_i}$ is a linear measure on $L_i$ for every $i \in I$. 
\end{defin}

Note that a linear measure on $\HH$ is in particular supported on a linear locus. 

\begin{prop}\label{prop: degree of an irreducible linear measure}
    An irreducible linear measure of dimension $d$ is homogeneous of degree $2d$. 
\end{prop}

\begin{proof}
    This is a consequence of the fact that if $\lambda$ is a Lebesgue measure on a complex linear space $V$, then for any $\C$-linear map $A : V \to V$, we have $A_*\lambda = |\det(A)|^{-2}\lambda$. We recall that the determinant of $A$ as a real linear transformation when $V$ is seen as a vector space over $\R$ is equal to $|\det(A)|^2$. 
\end{proof}

\begin{defin}
    A linear measure is said to be irreducible if its support is an irreducible linear locus. 
\end{defin}

\begin{defin}
    The {\em (complex) dimension} of an irreducible linear measure is the dimension of its support. 
\end{defin}

\subsection{Locally linear measures and linear measures on \texorpdfstring{$\HH$}{HH}}

This subsection is devoted to proving the following result which is a counterpart for measures of Proposition \ref{prop: locally linear locus is linear}:

\begin{prop}\label{prop: locally linear measure is linear}
    A locally linear measure on $\HH$ is a linear measure. 
\end{prop}

We will use the following lemma in the proof of Proposition \ref{prop: locally linear measure is linear}.

\begin{lem}\label{lem: restriction to linear submanifold}
    Let $m$ be a locally linear measure and let $L \subset \HHm$ be a linear submanifold which is contained in the support of $\pi^* m$ and such that $\pi^* m|_L$ is nonvanishing. Then $\pi^*m|_L$ is a linear measure supported on $L$. 
\end{lem}
 
\begin{proof}
    Let $\lambda$ be a Lebesgue measure on $V_L$ and let $m_L$ be the corresponding linear measure on $L$. We first show that for every $\q\in \HHm$, there is an open set $\U$ containing $\q$ and a constant $C> 0$ such that $\pi^* m|_{L \cap \U} = C \cdot (m_L)|_\U$. Indeed, let $(U,\Delta,\varphi)$ be a period chart around $\pi(\q)$ and let $\U$ be the neighborhood over $\varphi(U)$ such that $\varphi \circ \dev|_\U = \pi|_\U$. Since  $m$ is locally linear, there are linear subspaces $V_1,\dots,V_\ell \subset H^1(S,\Sigma;\R^2)$ such that $\varphi^*m$ is supported on $\bigcup_i V_i$ and $\varphi^*m|_{V_i}$ is the restriction of a Lebesgue measure $\lambda_i$ on $V_i$ for $i=1\dots \ell$. Without loss of generality, we can assume that the $V_i$ are pairwise distinct. We have by construction $(\dev|_\U)_*\big( (\pi|_\U)^*m\big) = \varphi^*m$. The assumption made on $L$ implies that $V_L \subset \cup_{i=1}^\ell V_i$. We have $\lambda_i|_{V_L} = 0$ unless $V_L=V_i$. Therefore, uniqueness of Lebesgue measure up to a multiplicative constant implies that there is a constant $C> 0$ such that $(\varphi^*m)|_{V_L} = C \cdot \lambda$. This implies that $(\pi^* m)|_{L \cap \U} = C \cdot (m_L)|_\U$.
       
    We now conclude by connectedness as follows. Cover $L$ by open sets $(\U_i)_{i \in I}$ such that for any $i \in I$, there is $C_i > 0$ satisfying $\pi^*m|_{L \cap \U_i} = C_i \cdot (m_L)|_{\U_i}$. Notice that $C_i = C_j$ whenever $\U_i \cap \U_j \neq \emptyset$. Define a function $f : L \to \R^+_*$ such that for any $\q \in U_i$, $f(\q) = C_i$. The function $f$ is locally constant and thus it is constant because $L$ is connected. This proves that $\pi^* m$ and $m_L$ coincide up to a constant.   
\end{proof}

\begin{proof}[Proof of Proposition \ref{prop: locally linear measure is linear}]
    Let $m$ be a locally linear measure on $\HH$. The support of a locally linear measure is a locally linear locus, and thus a linear locus by Proposition \ref{prop: locally linear locus is linear}. We denote it by $\M$. Let $\mathfrak{L}$ be a locally finite collection of linear submanifolds of $\HHm$ whose union is $\pi^{-1}(\M)$. Let $\mathfrak{L}_0$ be the sub-collection of linear submanifolds $L \subset \pi^{-1}\M$ such that $\pi^* m|_L$ is non-zero linear measure supported on $L$. We thus see that $\mathrm{supp} \ \pi^* m = \bigcup_{L\in \mathfrak{L}_0} L$ and for any $L \in \mathfrak{L}_0$, we have that $\pi^* m|_L$ is a non-zero linear measure supported on $L$ according to Lemma \ref{lem: restriction to linear submanifold}. 
\end{proof}

\subsection{Orbit closures and invariant measures described by linear loci and linear measures}

In this subsection, we use our notion of linear loci and linear measures to describe orbit closures and invariant measures in $\HH$.

\begin{thm}\label{thm: orbit is linear}
    The $\GL$-orbit closure of any $q \in \HH$ is an irreducible linear locus. 
\end{thm}

We first state and prove an intermediate lemma. 

\begin{lem}\label{lem: image of loc linear set}
    Let $\Gamma \subset \Mod(S,\Sigma)$ be a finite index, torsion free normal subgroup. The image of a  locally linear set in $\HH_\Gamma$ by the map $p_\Gamma : \HH_\Gamma \to \HH$ is a locally linear set in $\HH$.  
\end{lem}

\begin{proof}
    Let $\M_\Gamma$ be a locally linear set in $\HH_\Gamma$, let $\M = p_\Gamma(\M_\Gamma)$ and let $\Gamma_0 = \Gamma \backslash \Mod(S,\Sigma)$ (recall that $\Gamma$ is a normal subgroup). Up to replacing $\M_\Gamma$ by the union of its translates by $\Gamma_0$, which is again a locally linear set, we may assume that $p_\Gamma^{-1}(\M) = \M_\Gamma$. Let $(U,\Delta,\varphi)$ be a period chart for $\HH$. This chart can be lifted to a period chart for $\HH_\Gamma$, \textit{i.e} there is a chart $(U,\Delta_\Gamma,\varphi_\Gamma)$ for $\HH_\Gamma$ with the same domain $U$ such that $\varphi = p_\Gamma \circ \varphi_\Gamma$.  Notice that $\Delta_\Gamma$ is trivial since $\Gamma$ is torsion free. By the definition of a locally linear set, there are finitely many linear subspaces $V_1,\dots,V_\ell \subset H^1(S,\Sigma;\R^2)$ such that
    
    $$
        \varphi_\Gamma^{-1}(\M_\Gamma) = U \cap \bigcup_{i=1}^\ell V_i. 
    $$
    We thus compute that 
    
    $$
        \varphi^{-1}(\M) = \varphi_\Gamma^{-1}(\M_\Gamma) = U\cap \bigcup_{i=1}^\ell  V_i. 
    $$
    Since $(U,\Delta,\varphi)$ was arbitrary, this concludes our proof that $\M$ is a locally linear set. 
\end{proof}

\begin{proof}[Proof of Theorem \ref{thm: orbit is linear}]
    Let $\M$ be the orbit closure of $q$. Let $\Gamma \subset \Mod(S,\Sigma)$ be a finite index torsion free normal subgroup. Let $\Gamma_0$ be the finite group $\Gamma \backslash \Mod(S,\Sigma)$ and let $p_\Gamma : \HH_\Gamma \to \HH$ be the canonical projection map. Let $\tilde q \in \HH_\Gamma$ be a pre-image of $q$ by $p_\Gamma$ and let $\M_\Gamma$ be its orbit closure. It is a connected locally linear set by Theorem \ref{thm: emm} and its image by $p_\Gamma$ coincides with $\M$. Lemma \ref{lem: image of loc linear set} implies that $\M$ is locally linear and it is connected because $\M_\Gamma$ is and $p_\Gamma$ is continuous. We deduce from this that $\M$ is a linear locus by Proposition \ref{prop: locally linear locus is linear}. 
    
    We now show that $\M$ is an irreducible linear locus. Suppose  $\M=\M_1 \cup \M_2$ and $q \in \M_1$. By Proposition \ref{prop: invariant by the action of GL}, the $\GL$-orbit closure of $q$ is contained in $\M_1$ and thus $\M = \M_1$. This proves that $\M$ is an irreducible linear locus.  
\end{proof}

\begin{cor}\label{cor: irreducible linear locus is an orbit closure}
   Every irreducible linear locus in $\HH$ is a $\GL$-orbit closure.   
\end{cor}

\begin{proof}
    Let $\M$ be an irreducible linear locus in $\HH$ and let $\M_1$ be the set of surfaces in $\M$ whose $\GL$-orbit closure is strictly contained in $\M$. By \cite[Proposition 1.3]{EM} (see also \cite[Corollary 1.9]{Wright_cylinders}), there are at most countably many $\GL$-orbit closures in $\HH$ and by Theorem \ref{thm: orbit is linear}, each orbit closure is an irreducible linear locus. It follows that $\M_0$ is a countable union of positive codimension irreducible linear loci. In particular, $\M_0$ has empty-interior in $\M$ by the Baire Category Theorem. 
\end{proof}

Putting together Proposition \ref{prop:decomposin linear loci}, Theorem \ref{thm: orbit is linear} and Corollary \ref{cor: irreducible linear locus is an orbit closure}, we obtain a proof of the following result which was promised in \cite[Theorem 2.10]{SSWY}. 

\begin{thm}\label{thm: orbit is invariant subvariety}
    Irreducible linear loci and $\GL$-orbit closures coincide. Furthermore, any linear locus is a finite union of irreducible linear loci.
\end{thm}

We emphasize that invariant subvarieties as defined in \cite{SSWY} are examples of linear loci, but the converse does not hold in general. 

\begin{thm}\label{thm: invariant measure is cone measure associated to linear measure}
    Let $\nu$ be an ergodic $\SL$-invariant probability measure on $\HH^1$. Then there is $d>0$ such that $\mu_{2d} * \nu$ is an irreducible linear measure of dimension $d$ on $\HH$. 
\end{thm}

We first state and prove and intermediate lemma: 

\begin{lem}\label{lem: if the pre-image of a measure is locally linear}
     Let $\Gamma \subset \Mod(S,\Sigma)$ be a finite index torsion free normal subgroup and let $m$ be a measure on $\HH$. If the pre-image $p_\Gamma^*m$ of the measure $m$ by the quotient map $p_\Gamma : \HH_\Gamma \to \HH$ is a locally linear measure then the measure $m$ is a locally linear measure. 
\end{lem}

\begin{proof}
    Let $m_\Gamma$ be the pre-image of the measure $m$ by $\Gamma$. Let $(U,\Delta,\varphi)$ be a chart in $\HH$. This chart can be lifted, \textit{i.e} there is a chart $(U,\Delta_\Gamma,\varphi_\Gamma)$ with same domain $U$ such that $\varphi = p_\Gamma \circ \varphi_\Gamma$. Notice that $\Delta_\Gamma$ is trivial since we assumed that $\Gamma$ is torsion free. To conclude that $m$ is locally linear, it remains to observe that $\varphi^*m = \varphi_\Gamma^*m_\Gamma$. 
\end{proof}

\begin{proof}[Proof of Theorem \ref{thm: invariant measure is cone measure associated to linear measure}]
    Let $\nu$ be an ergodic $\SL$-invariant measure. Let $\Gamma \subset \Mod(S,\Sigma)$ be a finite index torsion free normal subgroup, let $\Gamma_0$ be the finite group $\Gamma \backslash \Mod(S,\Sigma)$ and let $p_\Gamma : \HH_\Gamma \to \HH$ be the canonical projection map. The action of $\Gamma_0$ on $\HH_\Gamma$ gives the map $p_\Gamma$ the structure of a quotient map as in Subsection \ref{quotientmaps}. Let $\nu_\Gamma$ be the pre-image of $\nu$ by $p_\Gamma$. By Proposition \ref{prop: image of a measure in the finite case}, we have ${p_\Gamma}_*\nu_\Gamma = |\Gamma_0| \cdot \nu$. The map $p_\Gamma$ is $\GL$-equivariant so $\nu_\Gamma$ is also $\SL$-invariant. Let $\nu_0$ be an ergodic component of $|\Gamma_0|^{-1} \cdot \nu_\Gamma$. We have $(p_\Gamma)_*\nu_0 = \nu$ and:
    $$
        \nu_\Gamma = \sum_{\gamma \in \Gamma_0} \nu_0 \prescript{}{*}{\gamma},
    $$ 
    where $\prescript{}{*}{\gamma}$ denotes the pushforward by the right action of $\gamma$.  Theorem \ref{thm: em} implies that there is $\delta>0$ such that $\mu_{\delta} *\nu_0$ is locally linear. We have: 
    $$
        \forall \gamma \in \Gamma_0, \ \left(\mu_{\delta} *\nu_0\right) \prescript{}{*}{\gamma} = \mu_{\delta} * \left( \nu_0 \prescript{}{*}{\gamma}\right).
    $$
    The translate of a locally linear measure by an element of $\Gamma_0$ is also a locally linear measure, and the class of locally linear measures is closed under finite sum. It follows that $\mu_{\delta} * \nu_\Gamma$ is a locally linear measure and since it coincides with the pre-image of $\mu_{\delta} * \nu$ by the map $p_\Gamma$, we deduce from Lemma \ref{lem: if the pre-image of a measure is locally linear} that $\mu_{\delta} * \nu$ is a locally linear measure. Proposition \ref{prop: locally linear measure is linear} now implies that $\mu_{\delta} * \nu$ is a linear measure. Let $\M$ be its support. It is a linear locus. Let $d = \dim_\C(\M)$. A linear measure of complex dimension $d$ is homogeneous of degree $2d$ and $\mu_{\delta} * \nu$ is homogeneous of degree $\delta$ so $\delta = 2d$. Assume $\M = \M_1 \cup \M_2$ can be written as a union of linear subloci, and assume $\mu_{\delta} * \nu(\M_1)>0$. It follows from Proposition \ref{prop : disintegration with respect to area function} that $\mu_{\delta} * \nu$ is also ergodic as a homogeneous measure and thus $\mu_{\delta} * \nu(\M\sm\M_1)=0$. By Proposition \ref{prop: invariant by the action of GL} we know that $\M_1$ is closed. The definition of the support of a measure implies that $\M \subset \M_1$. This proves that $\mu_{2d} * \nu$ is an irreducible linear measure.  
\end{proof}

\section{Linear loci, linear measures and geometric orbifolds}\label{sect: Linear loci, linear measures and geometric orbifolds}

In this section we discuss geometric structures on manifolds and orbifolds. We introduce special families of affine orbifolds which we call complex linear (resp. equilinear) orbifolds with a real structure, or simply real linear (resp. equilinear) orbifolds. We will describe linear loci and linear measures from the perspective of geometric structures. We then use the results of Section \ref{sec: recasting} to describe orbit closures and invariant  measures in $\HH$ from the perspective of geometric structures. 

\subsection{\texorpdfstring{$(G,X)$}{(G,X)}-structures on manifolds and orbifolds}\label{sec: linear structures}

Let $X$ be a manifold and let $G$ be a group acting smoothly and effectively on $X$. A $(G,X)$-structure on a manifold $M$ is given by an atlas of charts for $M$ taking values in $X$ so that the change of coordinate functions are restrictions of maps in $G$, see \cite[Chap. 5]{Goldman} for more details. For simplicity, we may omit the atlas from our notation and simply write $N$ for the $(G,X)$-manifold $(N,\mathcal{A})$ and write charts for $(G,X)$-charts. The connection between ``geometric structures" and geometry is the following. If there is some type of geometric object such as a metric, a volume form or a foliation defined on $X$ which is invariant under the action of $G$ then we can use the charts of our atlas to pull this back from $X$ to $M$. 

A simple example of such a structure is given by taking $G$ to be the additive group $\R^2$ and taking $X$ to be the space $\R^2$ where $G$ acts on $X$ by translations. The torus is an example of a space with this type of structure. We see that an $(\R^2,\R^2)$-structure on the torus gives rise to, for example, a horizontal foliation, an area form and a metric. It might seem reasonable to call an $(\R^2,\R^2)$-structure a translation structure but note that a general translation surface does not have an $(\R^2,\R^2)$-structure due to the presence of singular points.

As a second example of a $(G,X)$-structure we can take $X$ to be a vector space $V$ and $G$ to be the group of affine automorphisms of $V$, $Aff(V)$, then an $(Aff(V),V)$ structure is an {\em affine structure}. We can build a simple example of a one dimensional affine manifold by taking the quotient space of the positive real numbers by the action of the group $\{2^n: n\in\Z\}$. A similar construction in $\R^n$ gives us an example of a {\em Hopf manifold} (See \cite{Goldman} \S6.2). Typical affine structures differ from the ``translation structures" described previously in that they need not have translation invariant metrics. 

We introduced the developing map for strata in subsection \ref{subsec: strata}. This is a special case of a general construction for $(G,X)$ structures.
If $M$ is a $(G,X)$ manifold and $\widetilde M$ is its universal cover there is a developing map $dev: M\to X$ which is a $(G,X)$ map (\cite{Goldman} Theorem 5.2.2). This map determines the $(G,X)$ structure on $M$ in that we can pull back the structure from $X$ to $M$. 
There is a homomorphism $h:\pi_1(M)\to G$ so that for $\gamma\in\pi_1(M)$ so that for $\gamma\in\pi_1(M)$  
\begin{equation}\label{eq: holonomy}
    \begin{tikzcd}
        \widetilde M \arrow[r, "dev"] \arrow[d, "\gamma"] & X \arrow[d, "h(\gamma)"] \\
        \widetilde M \arrow[r, "dev"] & X
    \end{tikzcd}
\end{equation}

In fact according to (\cite{Goldman} Exercise 5.2.3) a developing map exists for any cover $\widetilde M$ for which $\pi_1(\widetilde M)$ is contained in the kernel of $h$. Developing maps give efficient ways of defining affine structures as we saw in subsection \ref{subsec: strata} where we used a developing map to define an affine structure on $\HHm$. In the setting of strata, holonomy is often called monodromy (see for example \cite{Calderon-Salter}).


Geometric structures can be defined on orbifolds as well as manifolds.
Orbifolds are discussed in, for example, Chapter 13 of Thurston's lecture notes \cite{thurston} and \cite{Caramello_orbifolds}.
It is useful when considering orbifolds to use an atlas of ``inverse charts" which map in the opposite direction from usual manifold charts. For the reader's convenience and to fix notation, we review orbifold structures and $(G,X)$-structures on orbifolds below.

An {\em orbifold} structure on a second countable, paracompact space $N$ is given by an atlas, which is a set $\mathcal{A}$ of compatible charts. A chart is a triple $(U, \Delta, \varphi)$ consisting of an open set $U$ in a fixed vector space $V$, a finite group $\Delta$ acting linearly and effectively on the set $U$ and a $\Delta$-invariant map $\varphi: U \to N$, which induces a homeomorphism from $U/\Delta$ onto $\varphi(U)$. Compatibility between charts is slightly more involved than in the case of manifolds: let us recall first that a chart $\alpha = (U_\alpha,\Delta_\alpha,\varphi_\alpha)$ embeds in a chart $\beta = (U_\beta,\Delta_\beta,\varphi_\beta)$ if there is a smooth embedding $\iota : U_\alpha \to U_\beta$ such that $\varphi_\alpha =   \varphi_\beta\circ\iota$. We say that $\iota$ is an embedding from $\alpha$ to $\beta$. Two charts $\alpha,\beta \in \mathcal{A}$ are compatible if for any $x \in \varphi_\alpha(U_\alpha) \cap \varphi_\beta(U_\beta)$, there is a chart $\gamma$ that embeds in both $\alpha$ and $\beta$ and such that $x \in \varphi_\gamma(U_\gamma)$. An atlas $\mathcal{A}_1$ refines an atlas $\mathcal{A}_2$ if every charts in $\mathcal{A}_1$ embeds in some chart in $\mathcal{A}_2$. This is an ordering on atlases and every atlas is refined by a unique maximal atlas and from now on we only consider maximal atlases. The local groups $\Delta$ give rise to a local group $\Delta_x$, depending only on $x \in N$ and which is well-defined up to a conjugation. If we require that domain of charts are instead open subsets of $X$ and that embeddings between charts and finite group actions are given by restrictions of maps in $G$, then we get a $(G,X)$-orbifold.

\begin{rem} 
    The paracompactness assumption ensures the existence of partitions of unity. These will play a role in the construction of a Lebesgue measure on real equilinear orbifolds, see Proposition \ref{prop: lebesgue measure on unimodular orbifold}.
\end{rem}

A simple example of an orbifold is the ``pillowcase" obtained by gluing two squares together (see Example 13.1.4 in Thurston's course notes). This space has a  $(G,X)$-structure where $G=\C$ acts on $X = \C$ by translations. We can view the pillowcase as the quotient of $\C$ by the group generated by $\Z^2$ and $z\mapsto -z$. This group $G$ preserves the quadratic differential $dz^2$ on $\C$ so the pillowcase has a canonical quadratic differential on it.
More generally any quadratic differential on the two sphere with four simple poles and no zeros has such a geometric orbifold structure. A general quadratic differential on a Riemann surface is not given by an orbifold structure due to the presence of zeros. Many examples of $(G,X)$-orbifolds arise as quotients, this is the point of Proposition \ref{prop: quotient orbifold} below. We first need the following definition.

\begin{defin}
    A map $f : N_1 \to N_2$ between $(G,X)$-manifolds is {\em a $G$-map} if it is locally of the form $x \mapsto x \cdot g$ for some $g \in G$ when expressed in charts. 
\end{defin}

\begin{defin}
    A map $f : N_1 \to N_2$ between $(G,X)$-manifolds is {\em an equivalence of $G$-structures} if $f$ and $f^{-1}$ are $G$-maps. 
\end{defin}

\begin{prop}\label{prop: quotient orbifold} 
    Let $N$ be a $(G,X)$-manifold and $\Gamma$ a group which is acting effectively, properly discontinuously on $N$ by $G$-maps. Then there is a unique $(G,X)$-orbifold structure on $N/\Gamma$ such that the quotient map $N \to N / \Gamma$ is a $G$-map.
\end{prop}

The proof of Proposition \ref{prop: quotient orbifold} is left to the reader. The fact that the quotient of a $(G,X)$-manifold by a properly discontinuous and free action of $\Gamma$ is a $(G,X)$-manifold is Exercise 5.1.1 in \cite{Goldman}. The fact that the quotient of a manifold by a properly discontinuous action is an orbifold is contained in Proposition 1.5.1 of \cite{Caramello_orbifolds}. 

We can interpret the results of Section \ref{sec: prelims} as showing that the manifold $\HH_{\mathrm{m}}$ has a $(G,X)$-structure where we take $G$ to be the trivial group and $X$ to be the vector space $H^1(S,\Sigma;\R^2)$. In particular, $\HHm$ carries an affine structure. An element $\gamma$ of the group $\Mod(S,\Sigma)$ acts on the model space $H^1(S,\Sigma;\R^2)$ by $\gamma^*$ which is a linear map. Therefore, it follows from Proposition \ref{prop: quotient orbifold} that strata have a well-defined affine orbifold structures. Note that strata provide examples of affine structures without translation invariant metrics. More on this can be found in subsection 2.3 of \cite{CSW}. 

Maps between orbifolds can be characterized in terms of {\em local lifts}. Let $f : P \to Q$ be a map between orbifolds, let $p \in P$ and let $q = f(p) \in Q$. A local lift of $f$ is the data of a chart $(U_1,\Delta_1,\varphi_1)$ around $p$, a chart $(U_2,\Delta_2,\varphi_2)$ around $q$ and a map $\tilde f: U_1 \to U_2$ so that the following diagram is commutative

\begin{equation}\label{eq: local lift}
    \begin{tikzcd}
        U_1 \arrow[r,"\tilde f"] \arrow[d, "\varphi_1"' ] & U_2 \arrow[d, "\varphi_2"] \\
        P \arrow[r, "f"] & Q
    \end{tikzcd}
\end{equation} 

For simplicity, we will say that $\tilde f$ is a local lift of $f$ and omit the charts, but it should be remembered that they are part of the definition. We say that $f$ is smooth if it has local lifts around every points that are smooth maps. See Section 1.4 of \cite{Caramello_orbifolds} for more details. We say that $f$ is an immersion when it has local lifts that are immersions. 

We will also make use of maps from $(G,X)$-orbifolds to $(G',X')$-orbifolds. This will depend on the choice of a family $M(X,X')$ of maps from $X$ to $X'$ which is closed under precomposition by elements of $G$ and postcomposition by elements of $G'$. It is then natural to consider maps that have local lifts in $M(X,X')$. A basic example would be to define affine maps between affine orbifolds by taking $(X,G)=(\R^m,Aff(R^m))$, $(X',G')=(\R^n,Aff(R^n))$ and $M(X,X')$ to be the family of affine maps from $\R^m$ to $\R^n$. More examples will be described in the next section. 

\subsection{Real (equi)linear structures}\label{sec: real linear structures}

The affine structures on strata have a number of distinctive properties. We can capture this fact by observing that strata have $(G,V)$-structures for a certain group $G$ which is a proper subgroup of $Aff(V)$.

As described in Definition \ref{defin: real structure} a real structure on a complex vector space $V$ is given by the choice of a real linear subspace $W\subset V$ so that $W\cap iW=\{0\}$ and $V=W\oplus iW$. If $V$ is a vector space with a real structure given by $W$ we denote by $\mathrm{GL}_\R(V,W)$ the group of complex linear automorphisms of $V$ that preserve $W$. 

\begin{defin}
    Let $V$ be a complex vector space with a real structure given by $W$. A {\em real linear orbifold} modeled on $V$ is a $(\mathrm{GL}_\R(V,W),V)$-orbifold. 
\end{defin}

We have described strata as being modeled on the vector space $H^1(S, \Sigma; \R^2)$. We can identify this with the complex vector space $V=H^1(S, \Sigma; \C)$. If we do so then $W=H^1(S, \Sigma; \R)$ is a real subspace which gives $\HHm$ a real linear manifold structure. 
Since the action of $\Mod(S,\Sigma)$ on $H^1(S, \Sigma; \C)$ leaves $H^1(S, \Sigma; \R)$ invariant the group $\Mod(S,\Sigma)$ acts on $\HHm$ by $\mathrm{GL}_\R(V,W)$-maps. This turns $\HH$ into real linear orbifold according to Proposition \ref{prop: quotient orbifold}. We can think of real linear orbifolds as ``stratum-like orbifolds". For describing invariant measures, it is useful to have an additional piece of structure.
We denote by $\mathrm{SL}_\R^{\pm}(V,W)$ the subgroup $\mathrm{GL}_\R(V,W) \cap \mathrm{SL}^{\pm}(V)$, where $\mathrm{SL}^{\pm}(V)$ is the subgroup of linear automorphisms of determinant $\pm 1$.  
 
\begin{defin}
   Let $V$ be a complex vector space with a real structure given by $W$. A {\em  real equilinear orbifold} modeled on $V$ is a $(\mathrm{SL}_\R^{\pm}(V,W),V)$-orbifold. 
\end{defin}

Equilinear orbifolds are closely related to the notion of affine manifolds with parallel volume that appears in \cite{Goldman}. The real linear manifold $\HHm$ has an equilinear structure coming from the developing map. If $\gamma$ is an element of the group $\Mod(S,\Sigma)$, then the action of $\gamma^*$ acts on the model space $H^1(S,\Sigma;\R^2)$ preserves the lattice $H^1(S,\Sigma;\Z^2)\subset H^1(S,\Sigma;\R^2)$ and so has determinant $\pm 1$. This means that strata are also real equilinear orbifolds. 
More generally the spaces $\HH_\Gamma$ for  $\Gamma \subset \Mod(S,\Sigma)$ a subgroup of finite index, are also real equilinear orbifolds. 

A map between real linear (resp equilinear) orbifolds is said to be {\em real linear (resp equilinear)} map if it has local lifts that are complex linear maps which preserves the real structure (resp. volume preserving complex linear maps which preserves the real structure). 

We now turn our attention to a special class of measures on real linear orbifolds which will be useful to us. In that perspective, Proposition \ref{prop: lebesgue measure on unimodular orbifold} will explain the advantage of real equilinear orbifold structures over real linear ones.   

\begin{defin}[Lebesgue measure on a real linear orbifold]
    Let $N$ be a real linear orbifold modeled on a linear space $V$. A Lebesgue measure on $N$ is a Radon measure $\lambda_N$ such that for every real linear chart $(U,\Delta,\varphi)$, there is a Lebesgue measure $\lambda$ on $V$ verifying 
        \begin{align}\label{eq: lebesgue measure in charts}
            \lambda_N|_{\varphi(U)} = \varphi_\ast (\lambda|_U). 
        \end{align}
\end{defin}

Orbifolds are locally described via finite-to-one quotient maps and constructing measures on the co-domain of quotient maps is precisely the point of subsection \ref{quotientmaps}. Proposition \ref{prop: image of a measure in the finite case} explains that in the finite case, the push-forward and the image of a measure by a quotient map only differ by a constant, and the class of Lebesgue measures on $V$ is preserved by multiplication by a positive scalar. That explains why the push-forward construction is sufficient in \eqref{eq: lebesgue measure in charts}. The Masur-Veech measure on $\HH$ is an example of a Lebesgue measure on a real equilinear orbifold. See \cite{Athreya-Masur} for more on the Masur-Veech measure. 

\begin{prop}\label{prop: lebesgue measure on unimodular orbifold}
     A real equilinear orbifold $N$ always carries a Lebesgue measure. When $N$ is connected, two Lebesgue measures only differ by a multiplicative constant. 
\end{prop}

The proof of Proposition \ref{prop: lebesgue measure on unimodular orbifold} will use the following lemma. 

\begin{lem}\label{lem: an technical lemma about equality of measure}
    Let $N$ be a real equilinear orbifold modeled on $V$, let $\lambda$ be a Lebesgue measure on $V$ and let $\iota : (U_1,\Delta_1,\varphi_1) \xhookrightarrow{} (U_2,\Delta_2,\varphi_2)$ be a chart embedding. Then the measure $\frac{1}{|\Delta_1|}(\varphi_1)_*\lambda$ coincides with the restriction to $\varphi_1(U_1)$ of $\frac{1}{|\Delta_2|}(\varphi_2)_*\lambda$.
    
\end{lem}

\begin{proof}
    The embedding $\iota$ is associated with an injective group morphism $\rho : \Delta_1 \to \Delta_2$ such that for every $\delta \in \Delta_1$ and $x \in U_1$, we have $\iota(x \cdot \delta) = \iota(x) \cdot \rho(\delta)$. Let $\delta_1,\dots,\delta_n \in \Delta_2$ be a maximal family of non-equivalent elements modulo $\rho(\Delta_1)$. We have
    $$  
        \varphi_2^{-1}\big(\varphi_1(U_1)\big) =  \iota(U_1) \cdot \Delta_2= \sqcup_{i = 1}^n \iota(U_1) \cdot \delta_i,   
    $$
    which follows from the fact that the $\Delta$ act effectively. Using the fact that $\iota$ is the restriction of a map preserving the measure $\lambda$ together with the equality $\varphi_1 = \varphi_2 \circ \iota$, we obtain: 
    \begin{align*}
        \frac{|\Delta_2|}{|\Delta_1|}(\varphi_1)_*\lambda = \frac{|\Delta_2|}{|\Delta_1|}(\varphi_2)_*\lambda|_{\iota(U_1)} = \sum_{i=1}^n  (\varphi_2)_*\lambda|_{U_1 \cdot \delta_i} = \big((\varphi_2)_*\lambda \big)|_{\varphi_1(U_1)}.\\
    \end{align*} 
\end{proof}

\begin{proof}[Proof of Proposition \ref{prop: lebesgue measure on unimodular orbifold}]
    Let $\mathcal{A}$ be the atlas defining the real linear structure on $N$ and let $\lambda$ be a Lebesgue measure on $V$. Cover $N$ by a locally finite collection of charts $\big\{(U_i,\Delta_i,\varphi_i)\big\}_{i \in I}$ in $\mathcal{A}$ (we use paracompactness here) and let $\{\rho_i : N \to \R^+\}$ be a partition of unity associated to this cover. We define a Radon measure $\lambda_{N}$ on $N$ satisfying the following formula
    $$
        \forall f \in C_c(N), \ \lambda_{N}(f) = \sum_{i \in I} \frac{1}{|\Delta_i|} \int_{U_i} f\rho_i \ d(\varphi_i)_*\lambda,
    $$  
    in light of the Riesz Representation Theorem. In order to show that $\lambda_N$ is indeed a Lebesgue measure, it is enough to show that for every chart $(U,\Delta,\varphi)$ in $\mathcal{A}$ and every $i \in I$, the measures $\frac{1}{|\Delta|}\varphi_*\lambda$ and $\frac{1}{|\Delta_i|}(\varphi_i)_*\lambda$ coincide on $\varphi(U_i) \cap \varphi(U)$. Lemma \ref{lem: an technical lemma about equality of measure} and compatibility of the corresponding charts imply that they coincide locally on $\varphi(U_i) \cap \varphi(U)$, which is enough to conclude. 
    
    Let now $\lambda_N'$ be another Lebesgue measure on $N$. Let $\alpha = (U_\alpha,\Delta_\alpha,\varphi_\alpha) \in \mathcal{A}$ be any chart. The Lebesgue measure on $V$ is unique up to a multiplicative constant so there is a constant $C_\alpha > 0$ such that $\lambda_{N}'|_{U_\alpha}=C_\alpha \cdot \lambda_{N}|_{U_\alpha}$. Furthermore, if $\alpha_1, \alpha_2 \in \mathcal{A}$ are such that $U_{\alpha_1} \cap U_{\alpha_2} \neq \emptyset$, then $C_{\alpha_1} = C_{\alpha_2}$. This means that there is a function $f : N \to \R^+_*$ such that for every $x \in N$ we have $f(x) = C_\alpha$ where $\alpha$ is any chart that contains $x$. By construction, this map is locally constant and thus constant if $N$ is connected. This ensures in that case that $\lambda_N$ and $\lambda_N'$ differ by a multiplicative constant. 
\end{proof}

\subsection{Linear loci and real linear orbifolds}\label{sec: Linear loci and complex linear orbifolds with a real structure}

In this subsection, we fix an irreducible linear locus $\M$. It follows from Proposition \ref{prop: max lin manifold} that if $L$ is a maximal linear submanifold in $\pi^{-1}(\M)$, we have $\pi(L) = \M$ and $\pi^{-1}(\M)$, which is the collection of translates of $L$ by elements of $\Mod(S,\Sigma)$, is locally finite. We called such an $L$ a lift of $\M$. We are going to derive some geometric consequences from the existence of lifts.  In particular we will give a description of the geometrical properties of $\M$ from the perspective of geometric orbifolds using a lift of $\M$. 

Let $L$ be a lift of $\M$, let $\mathrm{Stab}(L) \subset \Mod(S,\Sigma)$ be the stabilizer of $L$ and let $\mathrm{Stab}_0(L) \subset \mathrm{Stab}(L)$ be the normal subgroup of elements fixing $L$ pointwise. We denote by $\Gamma_L$ the quotient of $\mathrm{Stab}(L)$ by $\mathrm{Stab}_0(L)$. The group $\Gamma_L$ acts effectively and properly discontinuously on $L$, we denote by $\N_L$ its quotient and let $\pi_L : L \to \N_L$ be the canonical projection. The restriction to $L$ of the map $\pi$ induces a map $f_L : \N_L \to \HH$ whose image is $\M$ and that makes the following diagram commute: 

\begin{equation}\label{eq: diagram defining f_l}
    \begin{tikzcd}
        L \arrow[r, hookrightarrow] \arrow[d, "\pi_L"' ] & \HHm \arrow[d, "\pi"] \\
        \N_L \arrow[r, "f_L"] & \HH
    \end{tikzcd}
\end{equation}
By definition, $L$ is a real-linear manifold and $\Gamma_L$ acts by real linear maps. We endow the space $\N_L$ with the real linear orbifold structure described in Proposition \ref{prop: quotient orbifold}. With respect to this structure, the map $f_L$ is a real linear immersion. 

The action of $\GL$ on $L$ is given locally in charts by real linear maps. This action commutes with the action of $\Gamma_L$ and thus descends to an action of $\GL$ by real linear maps on $\N_L$ in such a way that the map $f_L$ is $\GL$-equivariant. 

Let $L^{\times} \subset L$ be the set of points in $L$ that are also contained in a distinct translate of $L$ by an element of $\Mod(S,\Sigma)$. Let 

\begin{equation}\label{eqref: singularity locus}
    \N_{L}^0 = \pi_L(L \sm L^\times). 
\end{equation}

\begin{prop}\label{prop: generically injective}
    The subset $\N_L^0$ is an open dense subset of $\N_L$.  The restriction of the map $f_L : \N_L \to \HH$ to the set $\N_L^0$ is injective.  
\end{prop}

\begin{proof}
    We first show that $L^\times$ is a nowhere dense closed subset of $L$. The intersection of $L$ with a translate of $L$ by $\gamma \in \Mod(S,\Sigma) \sm \mathrm{Stab(L)}$ is a locally linear submanifold of positive codimension $L$. It is thus both closed and nowhere dense in $L$. Furthermore, since $L$ locally encounters only finitely many translates of itself, we find that $L^{\times}$ is indeed closed and nowhere dense. This implies that $\N_L^0$ is an open dense subset of $\N_L$. 
    
    We now wish to show that $f_L$ is injective when restricted to $\N_L^0$. Say that two distinct points $q_1, q_2\in\N_L$ are mapped by $f_L$ to the same point in $\HH$. Let $\q_1$ and $\q_2$ be points in $L$ mapping to $q_1$ and $q_2$ by $\pi_L$. Since $\pi =f_L \circ\pi_L$, it follows that $\q_2 =\q_1 \cdot \gamma$ for some $\gamma\in\Mod(S,\Sigma)$.
    If $\gamma\in \mathrm{Stab}(L)$ then $\q_1$ and $\q_2$ map to the same point $\N_L$ under $\pi_L$ but this contradicts our assumption that $q_1 \ne q_2$. So $\gamma \notin \mathrm{Stab}(L)$. Now $\q_2 \in L$ and 
    $\q_2 = \q_1 \cdot\gamma  \in L \cdot \gamma$ so $\q_2 \in L^\times$ and thus $q_2 \notin \N_L^0$.
\end{proof}

\begin{prop}\label{prop: the map is proper}
    The map $f_L : \N_L \to \HH$ is proper. 
\end{prop}

\begin{proof}
    Let $K$ be a compact set in $\HH$. Let $\boldsymbol{K}$ be a compact set in $\HH_{\mathrm{m}}$ such that $\pi(\boldsymbol{K}) = K$. The collection of translates of $L$ by $\Mod(S,\Sigma)$ is locally finite so there is a finite collection  $\gamma_1,\dots,\gamma_\ell \in \Mod(S,\Sigma)$  such that $\boldsymbol{K} \cap \pi^{-1}(\M) \subset \bigcup_{i=1}^\ell L\cdot \gamma_i$. Define $\boldsymbol{K}_L = \bigcup_{i=1}^{\ell}
     \boldsymbol{K}\cdot \gamma_i^{-1} \cap L$. This is a compact set and we claim that $f^{-1}_L(K) = \pi_L(\boldsymbol{K}_L)$. This will imply that $f_L^{-1}(K)$ is also a compact set. To prove the claim, let $x \in \N_L$ be a point that is mapped in $K$ by $f_L$. Let $\q \in L$ be such that $\pi_L(\q) =x$. We have $\pi = f_L \circ \pi_L$ so there is $\gamma \in \Mod(S,\Sigma)$ such that $\q \cdot \gamma \in \boldsymbol{K}$. This implies that there exists $i = 1,\dots,\ell$ such that $\q \cdot \gamma \gamma_i^{-1} \in \boldsymbol{K}_L$ and $\pi_L(\q \cdot \gamma \gamma_i^{-1}) = \pi_L(\q) = x$. This proves the first inclusion. The other inclusion is immediate. Since $K$ was arbitrary, this proves that $f_L$ is indeed proper. 
\end{proof}

Recall that we say a map is {\em generically injective} when it is injective on an open dense subset of its domain.
Putting together the results of this subsection, we have proved the following: 

\begin{prop}\label{prop: linear locus implies linear orbifold}
    The orbifold $\N_L = L /\Gamma_L$ is a connected real linear orbifold endowed with an action of $\GL$ by real linear maps and the map $f_L : \N_L \to \HH$ is a proper and generically injective real linear equivariant immersion whose image coincides with $\M$. 
\end{prop}

The elements of this construction can be interpreted in the language of $(G,X)$ structures. In particular they may prove useful in calculating monodromy (holonomy). The space $L$ is a holonomy cover of $\N_L$. The developing map from $L$ to $V$ is obtained by restricting the developing map $dev:\HHm\to H^1(S,\Sigma;\R^2)$ to $L$. The image of the holonomy homomorphism is the image of the group $\Gamma_L$ under the relative homological representation.

\subsection{Linear measures and real equilinear orbifolds}

Let $m$ be an irreducible linear measure. Let $\M$ be the support of $m$. It is an irreducible linear locus. Let $L$ be a lift of $\M$ and let $f_L : \N_L \to \HH$ be the map described in \eqref{eq: diagram defining f_l}. Let $c = \# \mathrm{Stab_0}(L)$. The measure $c^{-1} \cdot \pi^*m$ is supported on $\pi^{-1}(\M)$ its restriction to $L$ is a $\Gamma_L$-invariant linear measure supported on $L$. Let $\lambda_{\N_L}$ be the induced measure on $\N_L$ \textit{i.e} the measure whose pre-image by the quotient map $\pi_L : L \to \N_L$ is $c^{-1} \cdot \pi^*m|_L$. The measure $\lambda_{\N_L}$ is a Lebesgue measure, and it follows from Proposition \ref{prop: degree of an irreducible linear measure} that it is a $\GL$-homogeneous Lebesgue measure of degree $2 \dim_\C(\M)$.

\begin{prop}\label{prop: linear measure as lebesgue measures on unimodular linear orbifolds}
    The space $\N_L = L /\Gamma_L$ is a connected real equilinear orbifold endowed with an action of $\GL$ by real linear maps, the map $f_L : \N_L \to \HH$ is a proper and generically injective real linear equivariant immersion whose image coincides with $\M$ and the measure $\lambda_{\N_L}$ is an ergodic Lebesgue measure such that $(f_L)_*\lambda_{\N_L} = m$. 
\end{prop}

\begin{proof}
    We already know by Proposition \ref{prop: linear locus implies linear orbifold} that $\N_L$ is a connected real linear orbifold endowed with an action of $\GL$ by real linear maps and that the map $f_L$ is a proper and generically injective real linear equivariant immersion whose image is $\M$. To prove that $\N_L$ is actually a real equilinear orbifold, it is enough to show that the action induced by $\Gamma_L$ on $V_L$, the model space of $L$ in period coordinates, preserves a Lebesgue measure. This is implied by the fact $\Gamma_L$ preserves the linear measure $\pi^*m|_L$. 
    
    Next, we show the formula $(f_L)_*\lambda_{\N_L} = m$. To that end, it is enough to prove that the pre-image by $\pi$ of $m$ and $(f_L)_* \lambda_{N_L}$ coincide on restriction to $L$. Let $\varphi \in C_C(L)$. We compute: 
    
    \begin{align*}
        \pi^* \left( (f_L)_* \lambda_{\N_L} \right)(\varphi) &= \int_\HH \sum_{\q \in \pi^{-1}(q)} |\Mod(S,\Sigma)(\q)|\varphi(\q) d{f_L}_*\lambda_{\N_L}(q) \\
        &= \int_{\N_L} \sum_{\q \in \pi^{-1}(f_L(q))} |\Mod(S,\Sigma)(\q)|\varphi(\q) d\lambda_{\N_L}(q) \\
        &= \int_{\N_L^0} \sum_{\q \in \pi_L^{-1}(q)} |\Mod(S,\Sigma)(\q)|\varphi(\q) d\lambda_{\N_L}(q) \\
        &= \int_{\N_L^0} \sum_{\q \in \pi_L^{-1}(q)} |\mathrm{Stab}(L)(\q)|\varphi(\q) d\lambda_{\N_L}(q) \\
        &= c \cdot \int_{\N_L} \sum_{\q \in \pi_L^{-1}(q)} |\Gamma_L(\q)|\varphi(\q) d\lambda_{\N_L}(q) \\
        &= c \cdot \pi_L^* \lambda_{\N_L}(\varphi) = \pi^*m(\varphi). 
    \end{align*}
    The third equality uses the fact $\lambda_{\N_L}(\N_L \sm \N_L^0) = 0$ because $L^\times$ is locally modeled on a finite union of positive codimension linear submanifolds. It also uses that for every  $q \in \N_L^0$, we have 
    $$
    \pi^{-1}(f_L(q)) \cap L =  \pi_L^{-1}(q).
    $$ 
    The fourth equality uses the fact that for every $\q \in \pi^{-1}(\N_L^0)$ we have the equality $\Mod(S,\Sigma)(\q) = \mathrm{Stab}(L)(\q)$.
    
    We finally prove that $\lambda_{\N_L}$ is ergodic. Let $B \subset \N_L$ be a $\GL$-invariant Borel set with $\lambda_{\N_L}(B) > 0$. The set $\N_L \sm \N_L^0$ is $\GL$-invariant and this is a $\lambda_{\N_L}$-conull set so we can assume that $B \subset \N_L^0$ and therefore $B = f_L^{-1}(f_L(B))$. The map $f_L$ is equivariant and so $f_L(B)$ is a $\GL$-invariant set. Since $(f_L)_*\lambda_{\N_L} = m$, we have $m(f_L(B)) >0$. Ergodicity of $m$ implies that $f_L(B)$ is a conull set for $m$ and so $B$ is a conull set for $\lambda_{\N_L}$.    
\end{proof}

\subsection{Orbit closures and invariant measures via real (equi)linear orbifolds}

We are finally in position to prove the main results announced in the introduction. 

\begin{proof}[Proof of Theorem \ref{thm: linear structure on orbit closures}]
   Theorem \ref{thm: orbit is linear} implies that a $\GL$ orbit closure is an irreducible linear locus $\M$. Let $L$ be a maximal linear locus in $\pi^{-1}(\M)$ then Proposition \ref{prop: linear locus implies linear orbifold} implies $\N_L$ is a real linear orbifold and $f_L:\N_L\to\HH$ is a proper and generically injective real linear orbifold immersion whose image is $\M$. 
\end{proof}

The construction of the real-linear structure on an orbit closure depends on the choice of a lift. We will show that the real-linear equivalence class of this real-linear structure is independent of these choices.

\begin{prop}
    Given an orbit closure $\M$ the construction of the corresponding $\N$ depends on the choice of a maximal linear manifold $L$ in $\pi^{-1}(\M)$. If $L$ and $L'$ are two such maximal linear manifolds then the corresponding spaces $\N$ and $\N'$ are not the same. We observe however that the real linear orbifolds $\N$ and $\N'$ are real linearly equivalent.
\end{prop}

\begin{proof} 
Since $\M$ is irreducible $L$ and $L'$ correspond under $\Mod(S,\Sigma)$. So $L'=L \gamma$ for some $\gamma\in\Mod(S,\Sigma)$ thought of as acting on $\HHm$ (on the right). The map $\gamma$ is a real linear equivalence of $L$ and $L'$. Let $\Gamma$ denote $\mathrm{Stab}(L)/\mathrm{Stab_0}(L)$ and $\Gamma'$ denote $\mathrm{Stab}(L')/\mathrm{Stab_0}(L')$ then $\Gamma'=\gamma^{-1}\Gamma\gamma$. The map $\gamma$ induces a real linear equivalence between $\N=L/\Gamma$ and $\N'=L'/\Gamma'$.
\end{proof}

\begin{proof}[Proof of Theorem \ref{thm: linear structure on invariant measures}]
     Theorem \ref{thm: invariant measure is cone measure associated to linear measure} gives an integer $d>0$ such that $\mu_{2d} * \nu$ is an irreducible linear measure of dimension $d$. By Proposition \ref{prop: degree of an irreducible linear measure}, the measure $\mu_{2d} * \nu$ is a homogeneous measure of degree $2d$ and by Proposition \ref{prop : coning is the same}, the measure $\nu$ is the cone measure associated to the linear measure $2d \cdot \mu_{2d} * \nu$. The result follows from Proposition \ref{prop: linear measure as lebesgue measures on unimodular linear orbifolds} applied to $2d \cdot \mu_{2d} * \nu$. 
\end{proof}

\end{document}